\documentclass[10pt,a4paper,leqno]{article}

\usepackage{amsmath,amsthm,amssymb,dsfont}
\usepackage{url,color,appendix,eucal,comment}
\usepackage{titlesec}
\usepackage[dvipsnames, table]{xcolor}
\usepackage{multicol}
\usepackage[T1]{fontenc}
\usepackage[all]{xy}
\usepackage{lscape,rotating}
\usepackage{graphicx}
\usepackage[hypertexnames=false]{hyperref}
\usepackage{tikz}
\usepackage{caption}
\usepackage{float}
\usepackage[lgrmath=cmr,lowercase=false,uppercase=false]{mathgreeks}

\newcommand{\GL}{\mathrm{GL}}
\newcommand{\ps}{\par \smallskip}
\newcommand{\N}{\mathbb{N}}
\newcommand{\Z}{\mathbb{Z}}

\newcommand{\Q}{\mathbb{Q}}

\newcommand{\R}{\mathbb{R}}
\providecommand{\C}{\mathbb{C}}

\newcommand{\F}{{\mathbb{F}}}

\newenvironment{pf}
{\medskip\noindent {\it Proof.  }}
{\hfill\nobreak $\Box$ \par\bigbreak}

\theoremstyle{plain}
\newtheorem{thm}{Theorem}[section]
\newtheorem{lemma}[thm]{Lemma}

\newtheorem{cor}[thm]{Corollary}
\newtheorem{prop}[thm]{Proposition}

\newtheoremstyle{plaincolon}
{}{}{\itshape}{}{\bfseries}{\,\,:}{.5em}{}
\theoremstyle{plaincolon}
\newtheorem*{conjHW}{Conjecture ${\rm HW(}Z{\rm )}$}
\theoremstyle{plain}

\newtheorem{thmintro}{Theorem}
\newtheorem{questionintro}{Question}
\newcounter{introcounter}
\theoremstyle{definition}
\newtheorem{definition}[thm]{Definition}
\newtheorem{remark}[thm]{Remark}
\newtheorem{example}[thm]{Example}

\titleformat{\section}{\bf \normalsize}{\arabic{section}.}{1 em}{}
\titleformat{\subsection}
{}
{\bf \arabic{section}.\arabic{subsection}.}
{0.5 em}
{\normalsize}

\titleformat{\subsubsection}[runin]{}{\thesubsubsection}{0.5em}{\itshape }

\hypersetup{colorlinks=true,linkcolor=black,citecolor=black,urlcolor=Periwinkle,pdfborderstyle={/S/U/W 1}}

\makeatletter

\@addtoreset{equation}{section}
\makeatother

\begin{document}

\title{On smooth projective hypersurfaces over $\Z$}

\author{
Ga\"etan Chenevier\thanks{\texttt{gaetan.chenevier@math.cnrs.fr}, CNRS, \'Ecole Normale Sup\'erieure-PSL, D\'epartement de Math\'ematiques et Application, 45 rue d'Ulm, 75230 Paris Cedex, France. During this work, G. Chenevier was supported by the project ANR-25-CE40-4664 PPAL.}
}


\maketitle

\begin{abstract} 
Assuming the Hasse-Weil conjectures for smooth hypersurfaces over $\Z$, 
we show that if there exists a degree $d$ hypersurface in $\mathbb{P}^{n-1}_{\Z}$
 that is smooth over $\Z$, with $n,d \geq 3$, then we have $n(d-2)/d > 555.6408$. 
If $n$ is odd, this last constant may be replaced with $661.168$. 
The proof combines the Riemann-Weil explicit formula with a certain weighted
Gaussian comparison for the distribution of Hodge numbers of smooth complex hypersurfaces.
We finally discuss automorphic variants of the converse.
\end{abstract}	
{\small
\tableofcontents
}
\ps
\ps

\section{Introduction}

\subsection{An open question} 
\label{subsect:openquest}
Let $n,d\geq 2$ be integers and let $F \in \Z[x_1,\dots,x_n]$ 
be a homogeneous polynomial of degree $d\geq 2$. 
The zero locus $F=0$ in the projective space $\mathbb{P}^{n-1}_\Z$ 
defines a projective hypersurface over $\Z$ of degree $d$ 
that we will denote by ${\rm Z}_F$. In this paper, 
we are interested in the following: 
	
\begin{questionintro} 
\label{quest:main}
For which integers $n,d \geq 2$ does there exist 
a hypersurface of degree $d$ in $\mathbb{P}^{n-1}_{\Z}$ that is smooth over $\Z$?
\end{questionintro}

For $F$ as above, recall that ${\rm Z}_F$ is smooth over $\Z$ 
if, and only if, for each prime $p$, the reduction mod $p$ of $F$ 
is nonzero and defines 
a nonsingular hypersurface in $\mathbb{P}^{n-1}(\overline{\mathbb{F}}_p)$. 
It is equivalent to ask that the {\it discriminant } of $F$, 
denoted by $\Delta$ in \cite[\S 2]{poonen} or ${\rm disc}$ in \cite{demazure}, 
is a unit in $\Z$, {\it i.e.} equal to\footnote{So equivalently, 
we are asking for the existence of an integer solution to $\Delta^2 =1$, 
a polynomial equation in the coefficients 
of the general hypersurface of degree $d$ in $n$ variables. 
This approach seems fairly impracticable though.}  
$\pm 1$. Note that if ${\rm Z}_F$ is smooth over $\Z$ then 
so is ${\rm Z}_{F'}$ for any $F'$ in the ${\rm GL}_n(\Z)$-orbit of $F$.
There are two classical cases where 
the answer to Question \ref{quest:main} is known. \ps

The first is the case of quadrics, {\it i.e.} $d=2$. 
Smooth quadrics over $\Z$ do exist in arbitrary dimension: 
we may take $F= \sum_{i=1}^{r} x_{2i-1}x_{2i}$ for $n=2r$ even, 
and $F= \sum_{i=1}^{r} x_{2i-1}x_{2i}+x_n^2$ for $n=2r+1$ odd. 
Better, it is well-known that 
the smooth quadrics over $\Z$ in $\mathbb{P}^{n-1}_\Z$ are exactly 
those defined by the rank $n$ even lattices of 
determinant $\pm 1$ (case $n$ even) or $\pm 2$ (case $n$ odd), 
two quadrics being in the same ${\rm GL}_n(\Z)$-orbit if and only if 
the associated lattices are isometric up to sign. 
The classification of these lattices is 
a well-known open problem in the definite case 
(which only happens for $n \equiv -1,0,1 \bmod 8$), 
but there is a single isometry class 
for each other possible signature otherwise: 
see e.g. \cite[Chap. V]{serre} and \cite[\S 2.2 \& Appendix B.2]{CL}. \ps

The second is the case of binary forms $F$, {\it i.e.} $n=2$. 
In this degenerate case, 
${\rm Z}_F \subset \mathbb{P}^1_\Z$ is $0$-dimensional over $\Z$, 
so it is smooth if and only if it is finite \'etale over $\Z$. 
A celebrated result of Minkowski 
(any nontrivial number field has a ramified prime) 
implies that each connected component of ${\rm Z}_F$ must have degree $1$ over $\Z$, so
${\rm Z}_F$ is the union of $d$ disjoint $\Z$-points in $\mathbb{P}^1_\Z$. 
Using that these points are disjoint over $\F_p$ for each $p$, 
we obtain first $d \leq |\mathbb{P}^1(\Z/2)|=3$, 
and for $d=3$ 
the only possibility is then ${\rm Z}_F = \{0,1,\infty\}$ and $F=x_1x_2(x_1-x_2)$, up to ${\rm GL}_2(\Z)$. 
Thus, for binary forms ($n=2$), a hypersurface as in Question~\ref{quest:main} exists exactly for $d\leq 3$.
\ps

As far as we know, there is no other example of 
a smooth hypersurface over $\Z$ in the literature, 
making the following (presumably folklore) variant 
of the question above especially interesting:

\begin{questionintro} 
\label{quest:folklore}
Do there exist integers $n,d \geq 3$ and a hypersurface of degree $d$ in $\mathbb{P}^{n-1}_{\Z}$ 
that is smooth over $\Z$ ?
\end{questionintro}

There are a few known negative results in this direction. 
The most famous one, 
due to Abrashkin and Fontaine \cite{abrashkin,fontaine}, asserts that 
{\it there is no smooth projective curve over $\Z$ of genus $>0$} (and actually no nonzero abelian variety over $\Z$), 
so that a positive answer to Question~\ref{quest:folklore} would require $n>3$. In another direction,
the question of obtaining lower bounds for the Artin conductor of 
projective smooth algebraic varieties over $\Q$ 
was studied from the perspective of the Riemann-Weil explicit formula 
by Mestre \cite{mestre}, in the spirit of classical work of Odlyzko \cite{odlyzko}. 
In particular, Mestre notes that the results of Fontaine and Abrashkin 
also simply follow from the conjectural analytic properties of 
the $L$ function of the $H^1$ of abelian varieties over $\Q$. 
This method was further developed in \cite{miller}, 
and then more recently in \cite{CL,chduke,chetai}, 
in which the authors obtain conditional (or automorphic) classifications of pure $\Q$-motives with good reduction everywhere 
and sufficiently small weights. For instance, 
it would easily follow from \cite[Thm. F]{CL}, 
together with a basic knowledge of the Hodge numbers of hypersurfaces 
(see below), that {\it a positive answer to Question~\ref{quest:folklore} 
should require\footnote{Let us mention that these arguments do not apply to the case $(n,d)=(4,3)$, 
but see Prop.~\ref{prop:cubisurfaces} for this folklore special case.}  $n > 24$}. \ps

\subsection{The main statement}
\label{subsect:mainstatement}
In this paper, we do not answer Question~\ref{quest:folklore}, 
alas!, but instead we study Question~\ref{quest:main} more thoroughly
from the point of view of the analytic methods above, 
especially \cite{mestre,chduke}. 
In particular, our first main result will be conditional on the standard conjectures on 
the Hasse-Weil zeta functions of hypersurfaces over $\Q$, 
that we now recall. \ps

Let $Z \subset \mathbb{P}^{n-1}_\Z$ be 
a hypersurface that is smooth over $\Z$. 
We know since Grothendieck and Deligne that, for each prime $p$, 
there is a polynomial $Q_p(t) \in 1+t \Z[t]$, 
with complex reciprocal roots of absolute value $p^{\frac{n-2}{2}}$, satisfying
$$\det(  1 - t \,{\rm Frob}_p\, \,|\,  
\,H^{n-2}_{\rm et}( Z_{\overline{\mathbb{F}}_p}, \Q_\ell ) )\,
= \,Q_p (t)\, \, {\rm \,\,for \,\,all\,\,}
\ell \neq p.$$
The $L$-function ${\rm L}(s,Z)$ is the Dirichlet series defined for $s \in \C$ with ${\rm Re}\, s> 1+(n-2)/2$
by the absolutely convergent Euler product
$${\rm L}(s,Z) \,\,= \,\,\prod_p \, 1/Q_p(p^{-s}).$$
On the other hand, there is a standard recipe 
associating an archimedean ${\rm L}$-factor $\Gamma(s, Z)$ 
to the complex Hodge structure $H^{n-2}(Z(\C),\C)$ equipped with its natural involution
\cite[\S 3.3]{serre_facteurs_locaux}, 
that we will recall in Sect.~\ref{subsec:gammafactors}.
Then the (unitarily normalized, complete) $L$-function of 
the degree $n-2$ cohomology of $Z$ is defined by
\begin{equation}
\label{eq:defLZ} \Lambda(s, \,Z) \,=\, \Gamma(s+(n-2)/2,Z) \,{\rm L}(s+(n-2)/2,Z)
\end{equation}
It converges absolutely on the half-plane ${\rm Re}\, s> 1 $. \ps

We can now state the Hasse-Weil conjecture for $Z$ in the form we shall use it, 
including the holomorphy expectation \`a la Artin and Tate \cite{tate_conjecture}. 
In the case $n$ is even, the conjecture involves the {\it rank} ${\rm r}_Z$ of $Z$, 
defined as the dimension of the subspace of $H^{n-2}(Z(\C),\Q)$ 
generated by the classes of the codimension $(n-2)/2$ algebraic cycles on $Z/\Q$.
Note that we always have ${\rm r}_Z\geq 1$ because of the Lefschetz class. 
We set ${\rm r}_Z=0$ for odd $n$.
We finally denote by $\xi(s) = \pi^{-s/2}\Gamma(s/2)\zeta(s)$ the classical complete Riemann zeta function. \ps

\begin{conjHW} 
\label{conj:conjconjHW}
The complete $L$-function $\Lambda(s, Z)$ has a meromorphic continuation to $\C$, with a functional equation of the form\footnote{The epsilon factor $\pm$ is purely archimedean here, and could be made entirely explicit.} $\Lambda(1-s, Z) = \pm \Lambda(s, Z)$. Moreover, $\Lambda(s,Z) \xi(s)^{-{\rm r}_Z}$ is entire and bounded in any vertical strip.
\end{conjHW}

This conjecture would follow for instance from the {\it automorphy conjecture} (Langlands, Fontaine-Mazur), 
asserting the existence of an isobaric automorphic representation 
$\pi$ of $\GL_b$ over $\Q$, with $b$ the $n-2$-th Betti number of $Z(\C)$, satisfying 
$\Lambda(s,\pi)= \Lambda(s,Z)$ and containing exactly ${\rm r}_Z$ times the trivial representation of $\GL_1$.
The main result of this paper is then the following:

\begin{thmintro}
\label{thmintro:HW} 
Assume $n,d \geq 3$ and let $Z$ be a hypersurface of degree $d$ in $\mathbb{P}^{n-1}_{\Z}$ that is smooth over $\Z$ and satisfies 
${\rm HW(}Z{\rm )}$. Then we have $n \frac{d-2}{d} > 661.168$ if $n$ is odd, 
and $n \frac{d-2}{d} >  555.6408$ otherwise. 
\end{thmintro}

In particular, {\it for $3 \leq n \leq 555$, and even $3 \leq n\leq 661$ for $n$ odd, and any $d\geq 3$, there should be no hypersurface of degree $d$ in $\mathbb{P}^{n-1}_{\Z}$ that is smooth over $\Z$}, in sharp contrast with the $n=2$ or $d=2$ cases. 
These large dimensions explain perhaps why we have never seen such a polynomial in the literature! 
As another example, {\it there should be no cubic hypersurface smooth over $\Z$ of odd dimension $\leq 1981$} (case $d=3$). 
\ps

\subsection{Ideas of the proof and organization of the paper}
\label{subsect:ideasandorg}

\begin{figure}[H]
\centering
\includegraphics[height=0.30\textheight,keepaspectratio]{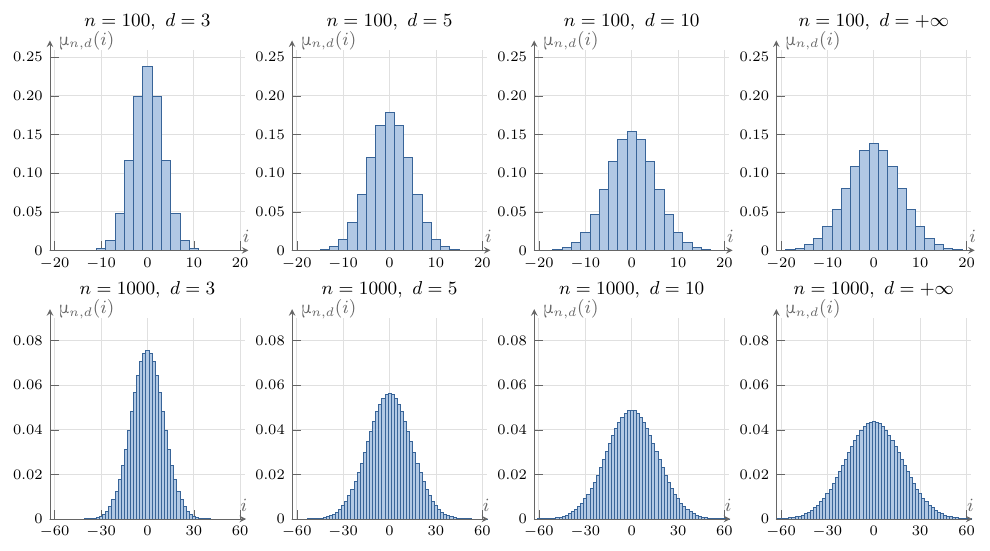}
\caption{The probability measure $\upmu_{n,d}$
for $n=100,1000$ and $d=3,5,10,+\infty$. }
\label{fig:hodge-walk-histograms}
\end{figure}

In Sect.~\ref{sect:hodgenumb}, we recall the Hodge numbers of 
a smooth complex hypersurface $Z$ of degree $d$ in $\mathbb{P}^{n-1}(\C)$, 
following Hirzebruch and Griffiths, but emphasizing a probabilistic point of view 
on the formula for these numbers as the probability law of a simple (conditioned) random walk. 
Precisely, for $i \in n+2\Z$ we denote by ${\rm h}_{n,d}(i)$ the middle degree, primitive, 
Hodge number of $Z$ with bidegree $((n-2+i)/2,(n-2-i)/2)$ (it does not depend on $Z$),
and by $\upbeta_{n,d}$ the primitive Betti number of $Z$ in degree $n-2$ (we have $\upbeta_{n,d} \neq 0$ for $n,d \geq 3$). Then 
\begin{equation}
\label{eq:defmundi}
\upmu_{n,d}(i) = {\rm h}_{n,d}(i)/\upbeta_{n,d}
\end{equation}
defines a probability measure on $\Z$ supported in $n+2\Z$. 
If $S_{n,d}$ denotes the random walk after $n$ steps
defined by the uniform law on $$\frac{1}{d}\{d-2,\,d-4,\,\dots,\,4-d,\,2-d\},$$
{\it then $\upmu_{n,d}$ is exactly the probability law of $S_{n,d}$
conditioned to lie in $n+2\Z$}. This interpretation, although immediate, does not seem to be widely known: 
the only reference to it we could find is\footnote{See also \cite[\S10.2]{lawrencevenkatesh} 
for a probabilistic study of the Eulerian distribution arising as $d \rightarrow + \infty$.} the recent work \cite[\S 3.1]{ji}.
The variance of $S_{n,d}$ is exactly 
\begin{equation}
\label{def:vnd}
{\rm v}_{n,d}=\frac{n(d-2)}{3d},
\end{equation} 
a quantity proportional to the one appearing in Theorem~\ref{thmintro:HW}. 
The exact variance of $\upmu_{n,d}$ is close to ${\rm v}_{n,d}$ for large $n$ and determined in \S\ref{subsect:Snd}. For a fixed $d\geq 3$ and $n \rightarrow \infty$, the central limit theorem shows thus that the law of $S_{n,d}/\sqrt{{\rm v}_{n,d}}$ 
weakly converges to the standard Gaussian measure $e^{-x^2/2}\,dx/\sqrt{2\pi}$ on $\R$,
and an elementary argument would show that the same actually holds for $\upmu_{n,d}$ rescaled by $\sqrt{{\rm v}_{n,d}}$, as illustrated in Fig.~\ref{fig:hodge-walk-histograms}.\footnote{We will not develop this aspect, as we do not need it, but we will prove instead an explicit ``weighted'' Gaussian comparison later in Theorem~\ref{thm:weightedgaussian}.}

In Sect.~\ref{sect:explform} we apply the {\it Explicit Formula} to $\Lambda(s,Z)$, following mostly \cite[\S I]{mestre}, and reduce the proof of Theorem~\ref{thmintro:HW} to that of an inequality that we now explain. An important role is played by an explicit numerical function that we call ${\rm J} : \R_{\geq 0} \rightarrow \R$, that is built from the family of test functions introduced by Odlyzko 
and the generalized digamma functions studied in \cite{mestre} and \cite{chduke}. We refer to \S\S\ref{sec:testdigamma} and~\ref{subsect:odlyzko} for the precise definitions.
For our purposes here, it is enough to say that ${\rm J}$ is a decreasing function, whose graph on $[0,50]$ is as in Fig.~\ref{fig:J} below, and
 with ${\rm J}(x)\,=\, {\rm log}\,4\pi\, -\, {\rm log}\,x\, +\, {\rm o}(1)$ for $x \rightarrow \infty$.
\begin{figure}[htbp]
\centering
\includegraphics[width=0.6\textwidth]{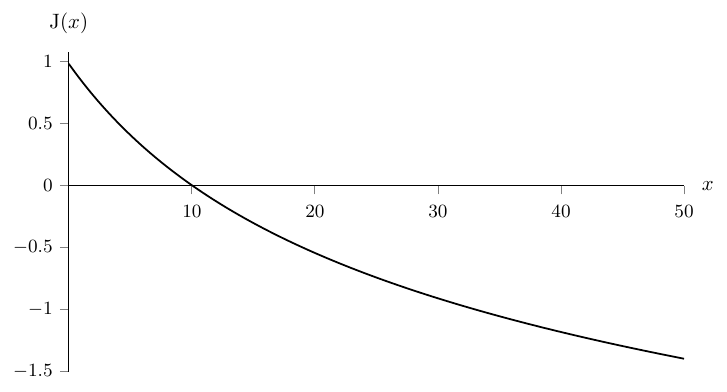}
\caption{The function $\mathrm{J}$ on $[0,50]$.}
\label{fig:J}
\end{figure}
The key object is then the expected value of ${\rm J} \circ |\cdot|$ with respect to $\upmu_{n,d}$, namely
\begin{equation}
\label{eq:defEnd}
{\rm E}_{n,d} \,=\, \sum_{i \in \Z}\, \upmu_{n,d}(i)\,{\rm J}(|i|) .
\end{equation}
If a $Z$ as in Theorem~\ref{thmintro:HW} exists, 
the Explicit Formula shows ${\rm E}_{n,d} \leq 0$ for $n$ odd, and ${\rm E}_{n,d} \leq \upmu_{n,d}(0) {\rm J}(0)$ for $n$ even (see Prop.~\ref{prop:maininequEF}).
The difference for $n$ even is that the classes in $H^{n-2}(Z,\Q)$ of the algebraic cycles of $Z$ of codimension $(n-2)/2$ and defined over $\Q$ contribute to poles of $\Lambda(s,Z)$, weakening the natural inequality given by the Explicit Formula (but see Remark~\ref{rem:algcaseneven}). \ps

It is important to remark that {\it for any given integer $d$, for instance $d=3$, Theorem~\ref{thmintro:HW} can now be proved\footnote{Strictly speaking the special case $(n,d)=(4,3)$ requires a specific treatment: Prop.~\ref{prop:cubisurfaces}.} by an explicit numerical check of the inequality ${\rm E}_{n,d}>0$ or ${\rm E}_{n,d} > \upmu_{n,d}(0) {\rm J}(0)$ for finitely many integers $n$.} This is how we actually first proved it for small $d$ and conjectured its general statement (with close constants). However, this numerical approach is not especially enlightening, and cannot establish the full theorem, since for any $n\leq 555$ {\it every} integer $d$ satisfies $n(d-2)/d \leq 555$. \ps

	The idea of our proof of Theorem~\ref{thmintro:HW} is as follows. Let $\uprho$ be the unique zero of ${\rm J}$. We have ${\rm J}(x) >0 \iff x<\uprho$ and a numerical computation shows
\begin{equation}
\label{eq:defrho}\uprho = 10.039802\ldots\, .
\end{equation}
It follows {\it e.g.} that ${\rm E}_{n,d}$ is trivially $>0$ for an odd $n<11$, proving Theorem~\ref{thmintro:HW} in this special case (a well-known fact). 
In general, the main idea is that, assuming that $\upmu_{n,d}$ approaches a normal distribution as above, and given the shape of the graph of ${\rm J}$, the positive values of ${\rm J}(|i|)$ for $|i| \leq 10$ will continue to dominate over the remaining negative values in ${\rm E}_{n,d}$ as long as ${\rm v}_{n,d}$ is not too large, which is exactly 
what is quantified in the statement of Theorem~\ref{thmintro:HW}. \ps

	In Sect.~\ref{sect:jensen}, we do not immediately follow this circle of ideas, but start by exploiting the straightforward lower bound ${\rm E}_{n,d} \geq {\rm J}(\sqrt{{\rm Var}\, \upmu_{n,d}})$, which follows from the Jensen inequality and the convexity of ${\rm J}$. It leads to a weaker but cheap version of Theorem~\ref{thmintro:HW}, with the already quite large constants\footnote{It is not a coincidence that $300.3928$ is close to $3 \uprho^2 -2$ !} $300.3928$ and $277.2971$ in place of $661.168$ and $555.6408$ respectively (Corollary~\ref{cor:weakmainthm}). The defect of this short argument is that it only uses the variance of $\upmu_{n,d}$ (and in the even case, an upper bound for $\upmu_{n,d}(0)$ as well), rather than the Gaussian shape of its distribution.
Nevertheless, it has a wider range of applications, as we shall see in Sect.~\ref{subsect:generaljensen}.
 \ps
     In the most technical Sect.~\ref{sect:weightedcomp}, we finally establish a version of the local limit theorem for $\upmu_{n,d}$ that is adapted to the weight ${\rm J}$ (and more generally to the family of weights ${\rm J}_H$ arising in the Explicit Formula). We assume $n$ odd now to simplify. The aim is to compare ${\rm E}_{n,d}$ with the Gaussian expected value
     $$\Theta_1(v) \,=\,  \sum_{i \in 2\Z+1} \frac{e^{-\frac{i^2}{2v}}}{{\rm m}_1(v)}\,\, {\rm J}(|i|),\,\,\,\,\,\,\,\, {\rm with}\, \, \,{\rm m}_1(v)=  \sum_{i \in 2\Z+1} e^{-\frac{i^2}{2v}},$$
at the point $v={\rm v}_{n,d}$. As we shall see, the map $v \mapsto \Theta_1(v)$ is decreasing on $\R_{>0}$, and admits a unique zero $\uprho_1$ satisfying $$3\uprho_1 = 661.3300065\ldots,$$
which is very close to the constant of Theorem~\ref{thmintro:HW}. The idea of course is that for $n$ large, ${\rm E}_{n,d}$ should be close enough to $\Theta_1({\rm v}_{n,d})$ so that they share the same sign, at least outside some small neighborhood of the sign change. \ps
	To make this idea effective, we thus 
provide upper and lower bounds on ${\rm E}_{n,d} - \Theta_1({\rm v}_{n,d})$ including explicit constants which are both sharp enough and uniform in $d$: see Theorem~\ref{thm:weightedgaussian}.  The arguments share a lot of similarities with the classical proofs of the local central limit theorems, going back to the pioneering work of de Moivre and Laplace (see e.g. \cite{feller, petrov}). But they also have some new features specific to the weight ${\rm J}$, which we discuss in Sect.~\ref{subsect:fourierreduction}. An important role there is played by a positive kernel reminiscent of that used in Lemma 2.3 and Proposition 2.2 in \cite{chduke}. This positivity, combined with  
a study of the sign of the first non-Gaussian term, explains {\it e.g.} why the value 
$3 \uprho_1$ above is slightly bigger than the $661.168$ in Theorem~\ref{thmintro:HW}.\ps
 The numerical application of our weighted local limit theorem to the proof of Theorem~\ref{thmintro:HW} is the subject of Sect.~\ref{sect:pfmainthm}. For $n$ odd, we actually prove that we have ${\rm E}_{n,d}>0 \iff n(d-2)/d<661.168$ (see Theorem~\ref{main:thm} (i)), showing that the bound $661.168$ is optimal, at least in this sense.\ps

\subsection{Perspectives, new questions and some answers}
\label{sect:perspectives}${}^{}$
	It is natural to ask whether the obstruction given by Theorem~\ref{thmintro:HW} 
to the existence of a degree $d$ hypersurface in $\mathbb{P}^{n-1}_{\Z}$ that is smooth over $\Z$ 
reflects an actual existence threshold, or if it is rather the byproduct of the Explicit Formula method that we use. 
The arguments above never really use the hypersurface itself, but merely the ${\rm L}$-function of its middle degree cohomology. 
This suggests considering the following variants of Question~\ref{quest:main}. 

\begin{questionintro} 
\label{quest:motiveQbar}
Assume $n,d\geq 3$. 
Does there exist a pure effective motive $M$ over $\Q$ of weight $n-2$, 
with coefficients in $\overline{\Q}$ and good reduction everywhere, 
whose Hodge numbers are those of the middle degree cohomology 
of a smooth degree $d$ hypersurface in $\mathbb{P}^{n-1}(\C)$ ?
\end{questionintro}

\begin{questionintro} 
\label{quest:motiveQ}
In the case Question~{\rm\ref{quest:motiveQbar}} has a positive answer, 
can such an $M$ be chosen with coefficients in $\Q$ 
and a nondegenerate pairing $M \otimes M \rightarrow \Q(2-n)$
that is alternating for $n$ odd and symmetric for $n$ even?
\end{questionintro}

We phrase these questions in terms of motives for convenience, 
but other formulations, either in terms of automorphic representations 
or in terms of ${\rm L}$-functions, would be equally interesting.
Of course, an affirmative answer to Question~\ref{quest:main} 
for $(n,d)$ implies an affirmative answer to Question~\ref{quest:motiveQ} for the same $(n,d)$. 
The converse however does not hold in general: 
it {\it e.g.} fails for $(n,d)=(4,3)$ since the Tate motive $\Q(-1)^{\oplus 7}$ satisfies all assumptions, 
whereas no cubic surface in $\mathbb{P}^3$ is smooth over $\Z$ 
(see Prop.~\ref{prop:cubisurfaces} and the footnote there). 
We leave open the question of which additional conditions could ensure the converse, 
as it seems to us to be a difficult problem 
and we have little to offer in this direction. \ps

Apart from the case $(n,d)=(4,3)$, our proof of Theorem~\ref{thmintro:HW} gives the exact same necessary 
inequality on $n(d-2)/d$ for the existence of an $M$ as in Question~\ref{quest:motiveQbar}, assuming\footnote{Our proof does not use the rationality of the coefficients of ${\rm L}(s,Z)$. Also,
the selfduality of ${\rm h}^{n-2}(Z)$ was only used to formulate the analytic assumptions in 
terms of a single ${\rm L}$-function rather than two.} 
the corresponding standard analytic conjectures on the complete normalized ${\rm L}$-functions $\Lambda(s,M)$ and $\Lambda(s,M^\vee)$ 
with respect to some embedding $\overline{\Q} \rightarrow \C$. 
Our aim now is to present a collection of positive results toward Questions~\ref{quest:motiveQbar} and~\ref{quest:motiveQ}.
We formulate them in automorphic terms, in order to avoid the conjectural passage from automorphic representations to motives. 
We refer {\it e.g.} to the preface in \cite{CL} and to \cite[\S 1.2]{chetai} for the background we need on level $1$ algebraic automorphic representations of $\GL_m$ and their relationship with Galois representations and motives. \ps

Assume $n,d \geq 3$. For simplicity, we also assume $n$ odd, 
although the case $n$ even could be treated by similar arguments. 
For all $f \geq 1$, we denote by $\Pi_{n,d,f}$ the set of isobaric irreducible automorphic representations of $\GL_m$ over $\Q$, with $m=\upbeta_{n,d}$, of the form $\pi = \pi_1 \boxplus \pi_2 \boxplus \cdots \boxplus \pi_s$,
for some $s\geq 1$, such that:\ps
-- for all $1 \leq j \leq s$, $\pi_j$ is a level $1$, algebraic regular, selfdual symplectic, cuspidal automorphic representation
of $\GL_{2m_j}$ over $\Q$, with $m_j \leq f$,\ps
-- for all $i \in \Z$, the multiplicity of the weight $i/2$ of $\pi_\infty$ is ${\rm h}_{n,d}(i)$.\ps
\noindent We finally denote by $\Pi_{n,d,f}^\Q \subset \Pi_{n,d,f}$ the subset of $\pi$ such that 
$\pi |\cdot |^{\frac{2-n}{2}}$ has rational Satake parameters.\footnote{By this we mean that for each prime $p$, the Satake parameter of $\pi_p |\cdot |^{(2-n)/2}$ has its characteristic polynomial with coefficients in $\Q$. For such a $\pi$, the Dirichlet series ${\rm L}(s-\frac{n-2}{2},\pi)$ has rational coefficients.} The following result is proved in Sect.~\ref{sect:persppf}. 

\begin{thmintro} 
\label{thm:autexist}
Assume $n$ is odd $\geq 3$.
\begin{itemize}
\item[(i)] We have $\Pi_{n,3,2} \neq \emptyset$ if, and only if, $n\geq 4893$.
\item[(ii)] If $n\geq 11139$, then we have $\Pi_{n,d,2} \neq \emptyset$ for every $d\geq 3$, and 
$\Pi_{n,d,2}^\Q \neq \emptyset$ for infinitely many $d$.
\item[(iii)] We have $\Pi_{2773,3,6} \neq \emptyset$.
\end{itemize}
\end{thmintro}

Concretely, assertion (i) shows that Question~\ref{quest:motiveQbar} should have a positive answer for $d=3$ and $n\geq 4893$ by simply considering motives that are direct sums of (suitable Tate twists of) symplectic motives of rank $2$ and $4$ attached respectively to modular eigenforms for ${\rm SL}_2(\Z)$ and Siegel modular eigenforms for ${\rm Sp}_4(\Z)$.\footnote{The construction of motives attached to cuspidal eigenforms for ${\rm SL}_2(\Z)$ 
  is a well-known result by Deligne and Scholl. The status of the analogous result for cuspidal eigenforms for ${\rm Sp}_4(\Z)$ is less clear to us, which explains our use of the conditional in the previous sentence.} 
  Using the same kind of motives, and for any fixed $n\geq 11139$, assertion (ii) predicts that Question~\ref{quest:motiveQbar} also has a positive answer for all $d \geq 3$, and that Question~\ref{quest:motiveQ} has a positive answer for infinitely many $d$. Last but not least, assertion (iii) provides an example with $n(d-2)/d=924+1/3$, which is not so far from our bound $661.168$. 
This suggests the following final amusing question:
  
 \begin{questionintro}
\label{quest:bestconstant}
{\it
What is the largest constant ${\rm C}_{\rm hyp}$ such that for all $n,d \geq 3$ with $n$ odd and $n(d-2)/d< {\rm C}_{\rm hyp}$, 
there is no motive $M$ as in Question~{\rm \ref{quest:motiveQbar}} ?
}
 \end{questionintro}  
  
The results of this paper show that ${\rm C}_{\rm hyp}$ should exist and satisfy 
\begin{equation}
\label{eq:intervallChyp}
661.168 \leq {\rm C}_{\rm hyp} \leq 924.333...,
\end{equation}
Although we do not know what ${\rm C}_{\rm hyp}$ should be, 
we certainly expect that the bounds above can be tightened a little 
using methods similar to the ones developed here. 
For instance, we already know that the lower bounds in Theorem~\ref{thmintro:HW} are not optimal, 
and can be improved under further analytic hypotheses. \ps

Indeed, the Explicit Formula method can also be applied to the Rankin-Selberg ${\rm L}$-function\footnote{An important difference is that these ${\rm L}$-functions may have a pole of high order at $s=1$.}  of the $H^{n-2}$ of the hypersurface $Z$, assuming a variant of ${\rm HW}(Z)$ for these ${\rm L}$-functions, or more strongly, the full automorphy of $H^{n-2}(Z)$. This technique has proved to give stronger constraints in related contexts ({\it e.g.} in \cite{CL,chetai}) and we checked indeed that combined with the classification results in \cite{CL}, it would allow us to increase the constants in Thm.~\ref{thmintro:HW} slightly in the odd case, and a bit more substantially in the even case (leading to a constant close to the one in the odd case). However, as this improvement is not so spectacular compared to the work needed to give the details, and does not really change the overall picture, we have chosen not to develop it here. In a similar direction, we mention that the use of GRH does not substantially increase our bounds either: see Remark~\ref{rem:grh}. \ps

All in all, it seems fair to say that Questions A and B remain entirely open in large dimension!
\bigskip\ps
{\sc Acknowledgments:} The author thanks Olivier Benoist, Djalil Chafa\"i and Olivier Ta\"ibi for useful discussions.\ps
\bigskip
{\sc AI declaration:}  This project is our first experiment with mathematical 
assistance from AI. In our original approach, we had discovered 
the statement of Theorem~\ref{thmintro:HW} 
using numerical experiments, proved it for small $d$, and identified a probabilistic 
approach based on a Gaussian comparison to explain this phenomenon. 
We took this opportunity to study with AI assistance the technical estimates 
needed to carry out this approach, building on Petrov's exposition 
of local central limit theorems \cite{petrov} and on our earlier work on explicit formulas \cite{chduke}. 
In the end, the arguments in Sect.~4 to 7 were developed in collaboration with OpenAI's GPT models,  
initially GPT-5.5, through several months of exchanges. We wrote the entire final text ourselves. 
The (straightforward) numerical verifications were run on PARI/GP \cite{pari}. The tables and figures were produced using code written by AI.
\ps

\section{The Hodge numbers of smooth hypersurfaces as a random walk}
\label{sect:hodgenumb}
\subsection{The Hodge numbers of smooth complex hypersurfaces}
\label{subsect:hodgehypersurfaces}
Let $n,d \geq 3$ be two integers.
Any homogeneous polynomial $F \in \C[x_1,\dots,x_n]$ of degree $d$ defines a hypersurface $Z \subset \mathbb{P}^{n-1}(\C)$, 
of dimension $n-2$.
Assuming $Z$ is smooth, then it is connected and
its Hodge numbers $h^{p,q}(Z)$ have been determined by Hirzebruch in \cite{hirzebruch}.
They can be derived as well from Griffiths' description of the primitive
 cohomology of $Z$ in terms of its Jacobian ring
\cite{griffiths}. That description shows that, for $p+q=n-2$, the {\it primitive, 
middle degree, Hodge number}\footnote{We set $\delta_P=1$ if the proposition $P$ holds, $\delta_P=0$ otherwise ({\it Kronecker symbol}).} 
\begin{equation}
\label{eq:primitivehodge}
{\rm h}_{n,d}^{p,q}=h^{p,q}(Z)-\delta_{p=q}
\end{equation}
is the coefficient of degree $(p+1)d$ in the polynomial 
$(t+t^2+\dots+t^{d-1})^n$, that is, {\it the number of ways to write $(p+1)d$ as a sum of $n$ integers between $1$ and $d-1$} (see {\it e.g.}  \cite[Cor. 17.5.4]{arapura}). Note that this number only depends on 
$d$, $n$ and $\{p,q\}$. 
We may easily recover from this formula the classical $n-2$-th Betti number of $Z$, 
namely $\dim H^{n-2}(Z,\C) = \upbeta_{n,d} +  \delta_{n \equiv 0 \bmod 2} $ with 
\begin{equation}
\label{eq:formbetand}
\upbeta_{n,d} = \frac{d-1}{d}\, \left( (d-1)^{n-1} - (-1)^{n-1} \right).
\end{equation} 
By our assumptions $n,d \geq 3$, we have $\upbeta_{n,d} > 0$. \ps

\subsection{A random walk interpretation}
\label{subsect:Snd}

We now explain the random walk interpretation of the Hodge numbers announced in Sect.~\ref{subsect:ideasandorg}. For our applications, it will be more natural to consider 
the {\it Hodge index difference} $i=p-q$. \ps

\begin{definition} 
\label{def:mund}
{\it For integers $n,d \geq 3$ and $-n<i<n$ with $i \equiv n \bmod 2$, 
we define the {\rm centered Hodge number} ${\rm h}_{n,d}(i)$ as ${\rm h}_{n,d}^{p,n-2-p}$ with $p:=\frac{n-2+i}{2}$, 
and set $$\upmu_{n,d}(i) = \frac{{\rm h}_{n,d}(i)}{\upbeta_{n,d}}.$$
We may and will view $\upmu_{n,d}$ as a discrete probability measure on $\Z$, whose support is included in $]-n,n[\cap (n+2\Z)$, and satisfying $\upmu_{n,d}(i) = \upmu_{n,d}(-i)$ for all $i \in \Z$.
 }
\end{definition}

Let $X_1^{(d)}, X_2^{(d)},\dots,X_n^{(d)}$ be independent 
random variables uniformly distributed in the discrete interval ${\rm I}_d=\{d-2,\, d-4,\,\dots,\,-d+2\}$ (``steps of size $2$''), and consider the associated random walk
$$S_{n,d}  = \frac{X_1^{(d)} + \dots + X_n^{(d)}}{d} \,\,\in \,\,\,]-n,n[.$$
We will condition this walk to the event $\mathcal{C}_{n,d} \subset {\rm I}_d^n$ 
defined by $S_{n,d} \,\in n+2\Z$.
As we shall see below, we have $\mathds{P}(\mathcal{C}_{n,d})>0$ for $n,d\geq 3$.

\begin{lemma} 
\label{lem:griffitsproba}
For $-n < i <n$ with $i \equiv n \bmod 2$,  we have 
$$\upmu_{n,d}(i) \,\, = \,\,\mathds{P}(\, S_{n,d} \,=\, i\, | \,\mathcal{C}_{n,d}).$$
\end{lemma}

\begin{proof} 
Indeed, the generating series of each $X_i^{(d)}$ is 
{\small
\begin{equation}
\label{eq:defpolPd}
{\rm P}_d(u)\,:= \,\frac{1}{d-1}( u^{d-2} + u^{d-4} + \cdots + u^{-d+4} + u^{-d+2} )  \,=\, \frac{ 1}{d-1} \frac{u^{d-1} - u^{-d+1}} {u - u^{-1}},
\end{equation}
}
so that of $dS_{n,d}$ is ${\rm P}_d(u)^n$. 
Note that for $p=\frac{n-2+i}{2} \in \Z$, we have $(p+1)d-n = (d-2)n/2+id/2$.
Using the identity
{\small
$$\frac{1}{(d-1)^n}\,(1+t+\dots+t^{d-2})^n\, t^{-\frac{(d-2)n}{2}} \,= {\rm P}_{d}(t^{1/2})^n,$$ 
}
\par \noindent Griffiths' formula reads ${\rm h}_{n,d}(i)=(d-1)^n{\mathds{P}}(S_{n,d} = i)$. 
By the obvious formula $\sum_i {\rm h}_{n,d}(i) = \upbeta_{n,d}$, we deduce
\footnote{This can also be checked directly, using {\it e.g.}  $\mathds{P}(\mathcal{C}_{n,d})=\frac{1}{2d} \sum_{\xi \in \upmu_{2d}} \xi^{-nd} {\rm P}_{d}(\xi)^{n}$.}
$\mathds{P}(\mathcal{C}_{n,d})=(d-1)^{-n} \upbeta_{n,d}>0$. 
\end{proof}

As the $X_i^{(d)}$ are uniform on ${\rm I}_d$ 
their average and variance are respectively $0$ and 
$\frac{1}{d-1} \sum_{i=1}^{d-1} (d-2i)^2 = \frac{d(d-2)}{3}$. 
It follows that we have $$\mathds{E}(S_{n,d})=0\, \, \, \,\,{\rm and}\, \, \, \,\,
{\rm Var}\, S_{n,d} = \frac{n(d-2)}{3d} = {\rm v}_{n,d}.$$
(Recall the notation ${\rm v}_{n,d}$ from Sect.~\ref{subsect:ideasandorg}.) 
The $\mathcal{C}_{n,d}$-conditioned average and variance of $S_{n,d}$ are as follows:

\begin{lemma} 
\label{lem:condvariance}
We have $\mathds{E}(S_{n,d}\, | \, \mathcal{C}_{n,d})=0$ and 
$${\rm Var}\, (S_{n,d}\, | \, \mathcal{C}_{n,d}) = 
(1 +  (-1)^{n+1} \frac{(n-1)(d-1)}{\upbeta_{n,d}})\,\,{\rm v}_{n,d} .$$
 \end{lemma}

\begin{proof} 
The symmetry of each $X_i^{(d)}$ and the equality $(n+2\Z)=-(n+2\Z)$ 
show $\mathds{E}(S_{n,d}\, | \, \mathcal{C}_{n,d})=0$. 
For the variance we have thus
$${\rm Var}\, (S_{n,d} \,| \, \mathcal{C}_{n,d}) \,=\, \mathds{E}( S_{n,d}^2\,|\,\mathcal{C}_{n,d})\,=\,
 \,\sum_{1 \leq i , j \leq n}  \,\mathds{E}(\frac{1}{d^2} X_i^{(d)} X_j^{(d)} \, |\, \mathcal{C}_{n,d}).$$
Write $X_i^{(d)}=d-2 Y_i$, with $Y_i$ uniform in $\{1,\dots,d-1\}$. 
In these variables, $\mathcal{C}_{n,d}$ coincides with $\sum_{i=1}^n Y_i \in d\Z$.
For $x \in \Z/d\Z$ and $m\geq 0$, consider 
$$\alpha_m(x)= \frac{1}{(d-1)^m}| \{ (x_1,\dots,x_m) \in (\Z/d\Z \smallsetminus \{0\})^m\, |\, \sum_{i=1}^m x_i = x\}|.$$
We have $\alpha_m(x)\,=\,\frac{1}{d(d-1)^m}\,\sum_{\xi \in \upmu_d} \xi^{-x} Q_m(\xi)$ with $Q_m=(t+t^2+\dots+t^{d-1})^m$, and $Q_m(\xi)=(-1)^m$ for $\xi \neq 1$, so we deduce:
\begin{equation}
\label{eq:formalphan}
\alpha_m(x) = \left\{ 
\begin{array}{cr} 
\frac{1}{d}(1+(-1)^m(d-1)^{1-m}) &   {\rm for} \, \, x=0, \\ \\
 \frac{1}{d}( 1-(-1)^m(d-1)^{-m})  = \alpha_{m+1}(0) & {\rm for}\,\, x \neq 0.
\end{array}
\right.
 \end{equation}
\noindent 
Of course, we have $\alpha_n(0)=\mathds{P}(\mathcal{C}_{n,d})$.
Moreover, as $\alpha_{n-1}(x)$ does not depend on $x$ for $x \neq 0$, 
the law of each coordinate of $(Y_1,\dots,Y_n)$ 
conditioned on $\mathcal{C}_{n,d}$ is uniform on $\{1,\dots,d-1\}$: we have for all $a \in \{1,\dots,d-1\}$
$$\mathds{P}(Y_i=a \,|\, \mathcal{C}_{n,d}) \,=\, \frac{1}{d-1} \alpha_{n-1}(-a) / \alpha_n(0)\, =\, \frac{1}{d-1}.$$ 
We have thus $\mathds{E}( \frac{1}{d^2} X_i^2 \,|\, \mathcal{C}_{n,d})=\mathds{E}( \frac{1}{d^2} X_i^2)$ for $i=1,\dots,n$, and 
$$\sum_{i=1}^n \mathds{E}( \frac{1}{d^2} X_i^2 \,|\, \mathcal{C}_{n,d}) =\sum_{i=1}^n \mathds{E}(  \frac{1}{d^2}  X_i^2 ) = {\rm Var}\, S_{n,d}.$$
We finally compute $\mathds{E}( X_i X_j \, |\, \mathcal{C}_{n,d})$ for $i \neq j$.
By definition, we have
$$ \mathds{E}( \frac{1}{d^2} \,X_i X_j \, |\, \mathcal{C}_{n,d}) \,=\, \frac{1}{d^2(d-1)^2} \sum_{1 \leq a,b \leq d-1} (d-2a)(d-2b) \frac{\alpha_{n-2}(-a-b)}{\alpha_n(0)}.$$
Note that the average of $(d-2a)(d-2b)$ in the sum above is zero, 
and that we have the identity $(d-2a)(d-2b)=-(d-2a)^2$ for $a+b=d$.
Using the identity
 $\alpha_{n-2}(x) = \alpha_{n-1}(0) + (\alpha_{n-2}(0)-\alpha_{n-1}(0)) \delta_{x=0}$, which follows from \eqref{eq:formalphan}, we finally find:
$$ \mathds{E}( \frac{1}{d^2}  X_i X_j \, |\, \mathcal{C}_{n,d}) = - \frac{1}{d^2(d-1)^2}  \frac{d(d-1)(d-2)}{3} \frac{\alpha_{n-2}(0)-\alpha_{n-1}(0)}{\alpha_n(0)}.$$
The statement follows then from $\alpha_{n-2}(0)-\alpha_{n-1}(0)=(-1)^n(d-1)^{2-n}$.
\end{proof}

\subsection{The $\Gamma$ factor of a smooth hypersurface defined over $\R$}
\label{subsec:gammafactors}
We now discuss the Archimedean $L$-factor $\Gamma(s,Z)$ attached to the degree $n-2$ cohomology
of a degree $d$ smooth hypersurface $Z$ in $\mathbb{P}^{n-1}_{\R}$, 
say for $n,d \geq 3$. This is a very special case of a general recipe going back 
at least to Serre in \cite{serre_facteurs_locaux}, which associates a $\Gamma$-factor to any 
finite dimensional graded complex vector space 
\begin{equation} \label{eq:CsurRHodge} H = \bigoplus_{(p,q) \in \Z^2} H^{p,q}\end{equation}
equipped with a $\C$-linear involution $\sigma$ satisfying $\sigma(H^{p,q})=H^{q,p}$ for all $(p,q) \in \Z^2$.
In our case, we have $H=H^{n-2}(Z(\C),\C)$ and \eqref{eq:CsurRHodge} is its Hodge decomposition, and 
$\sigma$ is the involution induced by the complex conjugation\footnote{This involution 
must not be mistaken with the semi-linear involution $x \mapsto \bar{x}$ on $H^{n-2}(Z(\C),\C)$ 
with fixed points $H^{n-2}(Z(\C),\R)$, and which also satisfies $\overline{H^{p,q}}=H^{q,p}$.} 
on $Z(\C)$ associated to the given real structure of $Z$. \ps

Recall the standard {\it real and complex $\Gamma$-factors}, defined as 
$\Gamma_\R(s) = \pi^{-s/2} \Gamma(s/2)$ and $\Gamma_\C(s) = 2(2\pi)^{-s} \Gamma(s)$,
where $\Gamma(s) = \int_0^{\infty} t^{s} e^{-t} \frac{dt}{t}$ is Euler's classical $\Gamma$-function. 
For $n$ odd, Serre's recipe simply takes the form\footnote{To compare with Serre's definition, just note $-(n-2)/2+{\rm min}(p,q)=-|p-q|/2$.}
\begin{equation}
\label{eq:gammanodd}
\Gamma(s+(n-2)/2,Z)\, = \prod_{\{p,q\} \,|\, p+q=n-2} \,\,\Gamma_\C(s+|p-q|/2)^{h^{p,q}(Z)},
\end{equation}
the product being over all pairs $\{p,q\} \subset \Z_{\geq 0}$ with $p+q=n-2$. 
Note that in terms of the Hodge index difference $i$ as above, 
we have $p=\frac{n-2+i}{2}$, $q=\frac{n-2-i}{2}$ and $p-q=i$. 
In this case $n$ odd, the real structure of $Z$ does not play any role. \ps

Assume now $n$ is even. The involution $\sigma$ induces a decomposition
$$H^{(n-2)/2,(n-2)/2} = H^+ \oplus H^-$$ where $(-1)^{(n-2)/2}\sigma$ acts as the identity on 
$H^+$ and as minus identity on $H^-$. Setting $h^{\pm} = \dim_\C H^\pm$, 
the general recipe for the $\Gamma$ factor of $H$ then reads:
\begin{equation}
\label{eq:gammaneven}
\Gamma(s+(n-2)/2,Z) = \Gamma_\R(s)^{h^+} \Gamma_\R(s+1)^{h^-} \prod \,\,\Gamma_\C(s+\frac{|p-q|}{2})^{h^{p,q}(Z)},
\end{equation}
the product being over all $\{p,q\} \subset \Z_{\geq 0}$ with $p+q=n-2$ and $p \neq q$. 
The quantities $h^+$ and $h^-$ are also uniquely determined by the relations
\begin{equation}
\label{eq:hpm}
\left\{\begin{array}{ccccc} h^+ + h^- & = & h^{(n-2)/2,(n-2)/2}(Z) & = & {\rm h}_{n,d}(0)+1 \\  h^+ - h^- & = & (-1)^{(n-2)/2}(\chi(Z(\R))-1)+1 \end{array}\right.,
\end{equation}
where $\chi(Z(\R))$ denotes the Euler characteristic of $Z(\R)$ (use e.g. \cite[Prop. I.2.1]{sihol}). By classical Hodge theory, if $Z$ is defined over $\Q$ we have the inequality
\begin{equation}
\label{eq:ineqrzhp} h^+ \geq {\rm r}_Z.
\end{equation}
\ps

\section{The Explicit Formula inequality}
\label{sect:explform}
Our aim in this section is to reduce the proof of Theorem~\ref{thmintro:HW} 
to two explicit inequalities. These inequalities will follow by studying the 
{\it Explicit Formula} for $\Lambda(s,Z)$, following mostly \cite[\S I]{mestre} 
(but see also \cite[\S 9.3.5]{CL} for several other examples). 
We start with 
a few preliminaries that will be important to state 
the archimedean contribution in this formula. \ps

\subsection{Test functions and generalized digamma functions}
\label{sec:testdigamma}
\begin{definition} 
\label{def:Hdig}
{\it Assume that $H : [0,+\infty[ \rightarrow \R$ is integrable and differentiable at $0$. 
 The associated {\it $H$-digamma} function is defined by the formula:
{\small
\begin{equation} \label{defpsih}
\begin{aligned}
\psi_H(s) & = \int_0^\infty \left( H(0)\frac{e^{-t}}{t} - H(t)\frac{e^{-st}}{1-e^{-t}} \right)\,{\rm d}t,\, \, \,\,\, {\rm Re}\, s \geq 0.
\end{aligned}
\end{equation}
}
\smallskip}
\end{definition}

This integral is absolutely convergent for all $s \in \C$ with ${\rm Re} \,s \geq 0$ by the assumptions on $H$,
 so $\psi_H(s)$ is a holomorphic function on ${\rm Re}\,s\,>0$. 
 Our terminology comes from the following examples (see also \cite[\S 2.1]{chduke}). \ps
 
\begin{example} {\it
\label{ex:Hdigamma} 
For $a>0$ and $H(t) = \frac{1-e^{-t}}{t}e^{-at}$ for $t>0$, 
a simple calculation shows $\psi_H(s)={\rm log}\, (a+s)$. 
Also, a classical formula due to Gauss shows $\psi_H(s)=\psi(a+s)$ for $H(t)=e^{-at}$, 
where $\psi(s)=\frac{\Gamma'(s)}{\Gamma(s)}$ is the classical {\it digamma} function. }
\end{example}
   
\begin{lemma} 
\label{lem:logcrois}
For any $H$ as in Definition~\ref{def:Hdig}, the map $s \mapsto \psi_H(s) - H(0) {\rm log}(1+s)$ is bounded on ${\rm Re} \,s \geq 0$ and converges to $0$ when ${\rm Re}\, s \rightarrow +\infty$. 
\end{lemma}

 \begin{pf} 
 In the case $H(0)=0$, we have $|\psi_H(s)| \leq \int_0^\infty \frac{|H(t)|}{1-e^{-t}} \,dt < \infty$ by the two assumptions on $H$, so the result follows from the dominated convergence theorem. In general, we reduce to this case by replacing $H(t)$ with 
 $G(t)=H(t)-H(0)\frac{1-e^{-t}}{t}e^{-t}$, since we have $\psi_{G}(s)=\psi_H(s) - H(0) {\rm log}(1+s)$ 
 by Example~\ref{ex:Hdigamma} and $G(0)=0$.  
 \end{pf}
 
	These $H$-digamma functions appear, following \cite{mestre}, 
when studying the contribution of $\Gamma$-factors to the Explicit Formulae with suitable {\it test functions}. 
For our purposes here, it will be enough to define a test function as 
{\it an even map $H : \R \rightarrow \R$ that is compactly supported and of class $\mathcal{C}^2$}. 
Such a function has a Fourier transform $\widehat{H}(\xi)=\int_\R H(x) e^{-2i\pi x \xi} dx$, 
an entire function, and an associated digamma function $\psi_H$.  \ps

Define the {\it group of $\Gamma$-factors} as the subgroup $\mathcal{G}$, 
of the multiplicative group of nonzero meromorphic functions on $\C$, 
generated by the functions of the form $a\,A^s$ or $\Gamma(as+b)$, with $a,A>0$ and $b\geq 0$. 
It contains thus all functions of the form $\Gamma_\R(as+b)$ and $\Gamma_\C(as+b)$ 
with $a>0$ and $b\geq 0$, by the duplication formula $\Gamma_\C(s)=\Gamma_\R(s)\Gamma_\R(s+1)$, 
hence all the $\Gamma$-factors defined in \S\ref{subsec:gammafactors} (after unitary normalization). 
A typical archimedean contribution in the explicit formula with test function $H$ will have the form 
\begin{equation}
\label{eq:defcalJH}
\mathcal{J}_H(\gamma)= - \int_\R \widehat{H}(t) \frac{\gamma'}{\gamma}(1/2+2i\pi t) dt
\end{equation}
for some $\gamma \in \mathcal{G}$. Such integrals have been studied by Mestre in \cite[\S I.2]{mestre}.
In particular they are absolutely convergent and real valued, and $\mathcal{J}_H : \mathcal{G} \rightarrow \R$ is (obviously) a group morphism. Moreover, Mestre gives {\it loc. cit.} the following alternative expressions, obtained by Fourier inversion:
 \begin{equation}
\label{eq:formJHdelta}
\left\{ \begin{array}{ccc} 
\mathcal{J}_H(\Gamma_\C(s+b/2)) & = & H(0) {\rm log} (2\pi) - \psi_H(\frac{1+b}{2}) \\
\mathcal{J}_H(\Gamma_\R(s)/\Gamma_\R(s+1)) & = & \int_0^\infty H(t) \frac{e^{-t/2}}{1+e^{-t}} dt
\end{array}
\right. .
\end{equation}

\par \noindent This leads us to introduce a last definition:

\begin{definition}
\label{def:JHdeltaH} {\it Assume $H : \R \rightarrow \R$ is integrable and differentiable at $0$. We set $\delta_H=\int_0^{\infty} H(t) \frac{e^{-t/2}}{1+e^{-t}} dt$, and define ${\rm J}_H : \R_{\geq 0} \rightarrow \R$ by the formula
\begin{equation}
\label{eq:defJH}
{\rm J}_H(x) =  H(0) {\rm log} (2\pi) - \psi_H(\frac{1+x}{2}).
\end{equation}
}
\end{definition}

\noindent The duplication formula $\Gamma_\C(s)=\Gamma_\R(s)\Gamma_\R(s+1)$ and Formula~\eqref{eq:formJHdelta}
also show
\begin{equation}
\label{eq:jh0delta}
2 \,\mathcal{J}_H(\Gamma_\R(s))= {\rm J}_H(0) + \delta_H\, \, \,{\rm  and}\, \, \, 2\, \mathcal{J}_H(\Gamma_\R(s+1))={\rm J}_H(0) - \delta_H.
\end{equation} 

\subsection{The Explicit Formula for $\Lambda(s,Z)$} ${}^{}$
\label{subsect:mainexplform}
Assume that $Z \subset \mathbb{P}^{n-1}_\Z$ is a degree $d$ hypersurface smooth over $\Z$, 
and that the conjecture ${\rm HW(}Z{\rm )}$ holds.
We fix a test function $H$ and set $\Phi_H(s) = \widehat{H}(\frac{1/2-s}{2i\pi})$ (an entire function of $s$).
We are going to state the explicit formula for $\Lambda(s):=\Lambda(s,Z)/\xi(s)^{{\rm r}_Z}$ and $H$. 
Roughly speaking, this is the identity that one obtains from integrating the $1$-form $\Phi_H(s)\,  {\rm dlog}\, \Lambda(s)$ 
along the boundary of the {\it critical strip} $0 \leq {\rm Re}\, s \leq 1$, 
and using both the residue theorem, the Euler product expansion of $\Lambda(s)$ and the functional equation $s \mapsto 1-s$.\ps

By definition, we may write $\Lambda(s)=\Lambda_\infty(s)\prod_p \Lambda_p(s)$ with 
$\Lambda_\infty(s)=\Gamma(s+(n-2)/2,Z)/\Gamma_\R(s)^{{\rm r}_Z}$ and $\Lambda_p(s) = (1-p^{-s})^{{\rm r}_Z}/Q_p(p^{-s-(n-2)/2})$. 
The running assumptions of \cite[\S I.2]{mestre} hold\footnote{In the notation of Mestre, we have 
$\Lambda_1(s)=\Lambda_2(s)=\Lambda(s)$ and $c=0$. 
The constants $A$ and $B$ appearing in Mestre's notation 
are purely coming from the $\Gamma$-factors in our presentations: 
the ``Artin conductor'' of $\Lambda$ is $1$ here by the formulation of ${\rm HW}(Z)$, 
but ultimately, since $Z$ is proper smooth over $\Z$. 
Note also that by standard properties of the $\ell$-adic cycle class map, 
the polynomial $(1-p^{(n-2)/2} t)^{{\rm r}_Z}$ divides $Q_p(t)$ in $\Z[t]$ for each prime $p$.}
by ${\rm HW}(Z)$, and the explicit formula takes the form:
\begin{equation}
\label{eq:explform}
\sum_{p,k}\, H(k \,{\rm log} \,p) \, a_{p^k}\,\frac{  {\rm log}\,p}{p^{k/2}}  \,+\, \mathcal{J}_H(\Lambda_\infty(s)) = - \frac{1}{2} \sum_{0 \leq \Re\, \rho \leq 1}\,\, \Phi_H(\rho)\,\, {\rm ord}_{s= \rho} \Lambda(s),
\end{equation}
where:\begin{itemize}
\item[(a)] (Contribution of finite primes) The (finite) sum on the left is over all primes $p$ and all integers $k\geq 1$, 
with $-\frac{\Lambda_p'}{\Lambda_p}(s) = \sum_{k\geq 1}  a_{p^k}\,\frac{{\rm log}\,p}{p^{ks}}$.\ps
\item[(b)] (Archimedean contribution) The term $\mathcal{J}_H(\Lambda_\infty(s))$ was defined in \eqref{eq:defcalJH}.\ps
\item[(c)] (Contribution of zeros) The (usually infinite) sum on the right is over all the zeros of $\Lambda(s)$ (since it has no pole by assumption), and has to be understood in the sense given in \cite[\S I.2]{mestre},\ps
\end{itemize}

\noindent This formula drastically simplifies under the following three assumptions on $H$:\ps
(A1) $H(x) \geq 0$ for all $x \in \R$,\ps
(A2) $H(x) = 0$ for $|x|\geq \log\,2$,\ps
(A3) ${\rm Re}\,\, \widehat{H}(\xi) \geq 0$ for\footnote{It is equivalent to ${\rm Re}\,\, \Phi_H(\xi) \geq 0$ for  $ 0 \leq {\rm Re}\, \xi \leq 1$.} all $\xi \in \C$ with $|{\rm Im}\, \xi|\leq \frac{1}{4\pi}$. \ps
\noindent Indeed, (the real part of) Formula \ref{eq:explform} then simply reads:
\begin{equation} 
\label{eq:JGammaisnegative}
\mathcal{J}_H(\Lambda_\infty(s)) \leq 0
\end{equation}

This is the analytic input we shall need. 
Recall the discrete $\upmu_{n,d}$ measure on $\Z$ defined in~\ref{def:mund}.
A major role will be played in the sequel by the following expected value, 
generalizing the quantity ${\rm E}_{n,d}$ in Formula~\eqref{eq:defEnd}. 

\begin{definition}
\label{def:EndH} {\it Assume $H : \R \rightarrow \R$ is integrable and differentiable at $0$ and let $n,d \geq 3$ be integers.
We set ${\rm E}_{n,d}^H = \sum_{i \in \Z} \upmu_{n,d}(i) {\rm J}_H(|i|)$. }
\end{definition}
In other words, we have ${\rm E}_{n,d}^H= \mathds{E}_{\upmu_{n,d}} ({\rm J}_H(|I|))$ where $I$ is any $\Z$-valued random variable with law $\upmu_{n,d}$. Recall that the support of $\upmu_{n,d}$ is finite and included in $[-(n-2),n-2] \cap (n+2\Z)$, in particular we have $\upmu_{n,d}(0)=0$ for $n$ odd.

\begin{prop} 
\label{prop:maininequEF}
Assume that there exists a degree $d$ hypersurface $Z \subset \mathbb{P}^{n-1}_\Z$ that is smooth over $\Z$ and satisfies 
{\rm (HW($Z$))}, with $n,d\geq 3$. Then for any test function $H$ satisfying {\rm (A1)}, {\rm (A2)} and {\rm (A3)} above we have ${\rm E}_{n,d}^H \leq \upmu_{n,d}(0) {\rm J}_H(0)$. 
\end{prop}

\begin{pf} 
Assume first $n$ is odd. 
Then using $\upmu_{n,d}(i)=\upmu_{n,d}(-i)$ 
together with Formulas~\eqref{eq:gammanodd},~\eqref{eq:formJHdelta} and~\eqref{eq:defJH}, 
we obtain $\frac{2}{\upbeta_{n,d}} \mathcal{J}_H(\Lambda_\infty(s)) = {\rm E}_{n,d}^H$, 
which is the statement by~\eqref{eq:JGammaisnegative}. 
Assume now $n$ is even. 
The same argument with Formula~\eqref{eq:gammaneven} instead 
shows that $\frac{2}{\upbeta_{n,d}} \mathcal{J}_H(\Lambda_\infty(s))$ coincides with
$$2\frac{h^+-{\rm r}_Z}{\upbeta_{n,d}} \mathcal{J}_H(\Gamma_\R(s)) + 2 \frac{h^-}{\upbeta_{n,d}} \mathcal{J}_H(\Gamma_\R(s+1)) + {\rm E}_{n,d}^H- \upmu_{n,d}(0) {\rm J}_H(0).$$
Note that if ${\rm J}_H(2) \leq 0$, 
then we have ${\rm J}_H(|i|)\leq 0$ for all even $i \neq 0$ 
as ${\rm J}_H$ is (clearly) decreasing, 
and the assertion trivially holds. 
So we may assume ${\rm J}_H(2) > 0$, 
but in this case observe we have 
${\rm J}_H(0)-{\rm J}_H(2) = \int_0^\infty H(t) e^{-t/2} dt  \geq \delta_H \geq 0$,
so ${\rm J}_H(0) \geq \delta_H$, 
and then both $\mathcal{J}_H(\Gamma_\R(s))$ and $\mathcal{J}_H(\Gamma_\R(s+1))$ 
are $\geq 0$ by Formulas~\ref{eq:jh0delta}. 
So the result follows again from \eqref{eq:JGammaisnegative}, 
together with $h^-\geq 0$ and \eqref{eq:ineqrzhp}.
\end{pf}

\begin{remark} \label{rem:algcaseneven} 
{\it
In the case $n$ odd, the conclusion of Prop.~\ref{prop:maininequEF} is ${\rm E}_{n,d}^H \leq 0$. 
The same bound can be obtained for $n$ even in the ``{\it balanced}'' case 
$h^+-h^-={\rm r}_Z=1$, or equivalently, ${\rm r}_Z=\chi(Z(\R))=1$ by~\eqref{eq:hpm}. 
Indeed, in this case the term
\begin{equation}
\label{eq:0partEF}
2 (h^+-{\rm r}_Z)\,\mathcal{J}_H(\Gamma_\R(s)) \,+ \,2 \,h^-\,\mathcal{J}_H(\Gamma_\R(s+1)) 
\end{equation}
coincides with $(h^++h^--1) {\rm J}_H(0) = {\rm h}_{n,d}(0) {\rm J}_H(0)$ by Formulas~\eqref{eq:jh0delta} and~\eqref{eq:hpm}. More generally, assuming\footnote{This is certainly satisfied for the interesting test functions as we shall see.} $\mathcal{J}_H(\Gamma_\R(s+1))\geq 0$, the term \eqref{eq:0partEF} is minimal (and actually zero) when $h^-=0$ and $h^+={\rm r}_Z$, so the bound stated in~Prop.\ref{prop:maininequEF} is the one obtained in this worst case scenario. We will not discuss here whether this case, or close ones, occur geometrically, but the point of this remark is to mention that the bound ${\rm E}_{n,d}^H \leq \upmu_{n,d}(0) {\rm J}_H(0)$ can be strengthened depending on extra information on $({\rm r}_Z,\chi(Z(\R)))$. }
\end{remark}

\subsection{Statement of the main theorem}
\label{subsect:odlyzko}
Since the pioneering work of Odlyzko \cite{odlyzko}, it is well-known that an especially efficient 
family of test functions satisfying (A1) and (A3) above is that of {\it Odlyzko functions}. 
They depend on a parameter $\lambda>0$ and are defined by 
\begin{equation}
\label{ed:defodl}
{\rm H}_\lambda(t) = \frac{{\rm g}(t/\lambda)}{\cosh(t/2)}
\end{equation}
where ${\rm g}: \R \rightarrow \R$ is defined by ${\rm g}(t)=(1-|t|)\cos(\pi t)+\pi^{-1}\sin(\pi|t|)$ for $|t| \leq 1$, and ${\rm g}(t)=0$ otherwise (see e.g. \cite[sect. 3]{poitou} and \cite[\S 9.3.17]{CL}). 
We have ${\rm H}_\lambda(0)=1$ and the support of ${\rm H}_\lambda$ is $[-\lambda,\lambda]$, 
so ${\rm H}_{\log\, 2}$ satisfies (A2).

\begin{definition} 
\label{def:J}
{\it In the special case $H={\rm H}_{{\rm log}\,2}$, we simply write ${\rm J}, \delta$ and ${\rm E}_{n,d}$ for ${\rm J}_H$, $\delta_H$ and ${\rm E}_{n,d}^H$ respectively.} 
\end{definition}

Our main numerical theorem will be the following.

\begin{thm}\label{main:thm} Assume $n,d\geq 3$. The following holds :
\begin{itemize}
\item[(i)] For $n$ odd, we have ${\rm E}_{n,d} >0\,\,\iff\,\,n\frac{d-2}{d} < 661.168$.\ps
\item[(ii)] For $n$ even with $(n,d)\neq (4,3)$ and $n\frac{d-2}{d} < 555.6408$ we have ${\rm E}_{n,d}>\upmu_{n,d}(0){\rm J}(0)$.
\end{itemize}
\end{thm}

The proof of this theorem will occupy most of what follows. 
It implies Theorem~\ref{thmintro:HW} by Prop.~\ref{prop:maininequEF}, 
except in the special case $(n,d)=(4,3)$ 
(for which the method above is empty as we have ${\rm r}_Z=7$ and $\Lambda(s,Z)=\xi(s)^7$). 
It is however a folklore result that this special case does not occur:

\begin{prop}
\label{prop:cubisurfaces}
 (Folklore) There is no smooth cubic surface over $\Z$. 
 \end{prop}

\begin{pf} 
A sketch of proof is given by A. Javanpeykar in 
\href{https://mathoverflow.net/questions/336618/good-reduction-of-rational-surfaces}{this MathOverflow post}. 
Here is a variant of this proof. 
If $S$ is a smooth cubic surface in $\mathbb{P}^3_{\Q}$, 
then we know that $H^2_{\rm et}(S_{\overline{\Q}},\Q_\ell(1))$ is 
$7$-dimensional and spanned by algebraic classes, 
so that its natural Galois action factors through a finite quotient. 
If $S$ admits furthermore a proper smooth model over $\Z$, 
we know that this action is both unramified outside $\ell$ and crystalline at $\ell$ 
(hence unramified as well since the image is finite), 
so it is trivial by Minkowski's theorem 
(showing actually $\Lambda(s,S)=\xi(s)^7$ as claimed above). 
The Grothendieck-Lefschetz formula then reads $|S(\mathbb{F}_p)|=1+7p+p^2$. 
But for $p=2$ this is $19$, which is greater than $|\mathbb{P}^3(\F_2)|=2^4-1=15$. 
So $S$ does not embed\footnote{As Javanpeykar points out {\it loc. cit.}, 
beware that there do exist smooth cubic surfaces over $\Q$ 
having a proper smooth model over $\Z$ 
(but this model is not a hypersurface in $\mathbb{P}^3$ over $\Z$).} 
in $\mathbb{P}^3_{\Z}$.
\end{pf}

\section{A cheap weaker result: the Jensen bound}
\label{sect:jensen}
The statements in Theorem~\ref{main:thm} concern the correlation between 
$\upmu_{n,d}$ and the explicit function ${\rm J} : \R_{\geq 0} \rightarrow \R$ defined in \ref{def:JHdeltaH} and \ref{def:J}. In this section, we will focus on some numerical properties of ${\rm J}$, but mostly ignore the 
specificities of the probability measure $\upmu_{n,d}$: all we shall need about it is its variance (from Lemma~\ref{lem:condvariance}), as well as an upper bound on $\upmu_{n,d}(0)$ for $n$ even.

\subsection{The function ${\rm J}$} We refer to \cite[\S 9.3.17]{CL} for some numerical methods to compute ${\rm J}$ up to arbitrary precision, 
and first recall a few observations made there. Here are a few values.

\begin{table}[htbp]
\centering
\small
\setlength{\tabcolsep}{2.5pt}
\begin{tabular}{c|*{10}{r}}
$i$ & 0 & 1 & 2 & 3 & 4 & 5 & 6 & 7 & 8 & 9 \\
$\mathrm{J}(i)$
& $0.982$ & $0.848$ & $0.725$ & $0.611$ & $0.506$
& $0.408$ & $0.316$ & $0.230$ & $0.150$ & $0.074$ \\
\hline
$i$ & 10 & 11 & 12 & 13 & 14 & 15 & 16 & 17 & 18 & 19 \\
$\mathrm{J}(i)$
& $0.003$ & $-0.065$ & $-0.129$ & $-0.189$ & $-0.247$
& $-0.302$ & $-0.355$ & $-0.405$ & $-0.453$ & $-0.499$
\end{tabular}
\caption{Values of $\mathrm{J}(i)$ rounded to three decimal places.}
\label{tab:J}
\end{table}

\begin{prop} The function ${\rm J}$ is smooth, strictly decreasing, and convex on $\R_{\geq 0}$, and it satisfies 
${\rm J}(x) = {\rm log}(4\pi) - {\rm log}(x) + {\rm o}(1)$ when $x \rightarrow \infty$.
\end{prop}

\begin{pf} For any test function $H$, the function $\psi_H(s)$ is entire on $\C$ and satisfies 
$\psi_H'(s) = \int_0^\infty  t H(t) e^{-st}/(1-e^{-t}) dt$ and $\psi_H''(s) = \int_0^\infty - t^2 H(t) e^{-st}/(1-e^{-t}) dt$, 
by the dominated convergence theorem. It follows that for $H\geq 0$, the map ${\rm J}_H(x)=H(0){\rm log}(2\pi)-\psi_H((1+x)/2)$ is smooth, decreasing (strictly if $H\neq 0$) and convex on $\R_{\geq 0}$. The last assertion follows from Lemma~\ref{lem:logcrois}.
\end{pf}

This shows that ${\rm J}$ has a unique zero on $\R_{\geq 0}$, that we denoted by $\uprho$ 
in \S\ref{subsect:ideasandorg}. 
A numerical computation using the formulas in \cite[\S 9.3.17]{CL} shows $\uprho = 10.0398024452...$.
\subsection{The Jensen bound}
\label{subsect:jensen}
Our aim now is to prove the following weak, but cheap, version of Theorem~\ref{main:thm}.

\begin{prop} 
\label{prop:inegjensen}
Assume $n,d \geq 3$ and $(n,d)\neq (4,3)$. \begin{itemize}
\item[(i)] If $n$ is odd and $\frac{n(d-2)}{d} < 300.3928$, then we have ${\rm E}_{n,d} >0$.\ps
\item[(ii)] If $n$ is even and $\frac{n(d-2)}{d}<277.2971$, then we have ${\rm E}_{n,d} > \upmu_{n,d}(0) {\rm J}(0)$.
\end{itemize}
\end{prop}

\noindent Using Propositions~\ref{prop:maininequEF} and \ref{prop:cubisurfaces}, we immediately deduce:

\begin{cor} 
\label{cor:weakmainthm}
The statement of Thm.~\ref{thmintro:HW} holds with $661.168$ and $555.6408$ replaced respectively with 
$300.3928$ and $277.2971$. 
\end{cor}

We now prove Proposition~\ref{prop:inegjensen}. Recall the notation ${\rm v}_{n,d}=\frac{n(d-2)}{3d}$.

\begin{pf} Let $I$ be a $\Z$-valued random variable with law $\upmu_{n,d}$. We have 
$$\mathds{E}(\,{\rm J}(|I|)\,) \geq {\rm J}(\,\mathds{E}(\,|I|\,)\,) \geq  {\rm J}(\,\sqrt{{\rm Var}\, I}\,).$$
Indeed, the first inequality follows from the convexity of ${\rm J}$ and the Jensen inequality, 
and the second from the Cauchy-Schwarz inequality $\mathds{E}(|I|) \leq \sqrt{ {\rm Var}\,I }$ as ${\rm J}$
is decreasing. As ${\rm J}$ is strictly decreasing again, we have ${\rm J}(\sqrt{{\rm Var}\, I})>0$ if, and only if, ${\rm Var}\, I<\uprho^2$. Recall the formula 
\begin{equation}
\label{eq:recallvar} {\rm Var}\, I \,=\, {\rm v}_{n,d}\,(1+\epsilon_{n,d})\, \, {\rm with}\,\, \epsilon_{n,d}=(-1)^{n+1}(n-1)(d-1)/\upbeta_{n,d}
\end{equation}
from Lemma~\ref{lem:condvariance}. For $n$ odd, we also have ${\rm Var}\, I \leq {\rm v}_{n,d}+2/3$ by Lemma~\ref{lem:majvarodd} below, so assertion (i) of Proposition~\ref{prop:inegjensen} follows from the numerical value
$$ 3\uprho^2-2 = 300.3928994...$$ 

We finally assume $n$ is even. The argument above, together with the bound ${\rm Var}\, I \,\leq\, {\rm v}_{n,d}$ from Lemma~\ref{lem:majvarodd}, still shows ${\rm E}_{n,d}>0$ for $n(d-2)/d < 3 \uprho^2$, but this is not what we need.
Assume furthermore $(n,d) \neq (4,3)$. By the following Lemma~\ref{lem:hnd2nonzero}, 
note that we have $\upmu_{n,d}(0) < 1$.
Let $\upmu^*_{n,d}$ be the conditional law of $I$ given $I \neq 0$, and let $I^*$ be a $\Z$-valued random variable with law $\upmu^*_{n,d}$. We simply have $\upmu^*_{n,d}(0)=0$ and $\upmu_{n,d}^*(i)=\upmu_{n,d}(i)/(1-\upmu_{n,d}(0))$ for $i \neq 0$. 
By the same arguments as above we get
$$\frac{{\rm E}_{n,d}-\upmu_{n,d}(0){\rm J}(0)}{1-\upmu_{n,d}(0)} = \mathds{E}(\,{\rm J}(|I^*|)\,) \geq {\rm J}(\,\mathds{E}(\,|I^*|\,)\,)  \geq  {\rm J}(\,\sqrt{{\rm Var}\,I^*}\,).$$
The variable $I^*$ is clearly centered with variance ${\rm Var}\,(\,I^*) = ({\rm Var}\, I)/(1-\upmu_{n,d}(0))$. 
It follows that if we have ${\rm Var}\, I < (1-\upmu_{n,d}(0)) \uprho^2$, then ${\rm E}_{n,d}-\upmu_{n,d}(0){\rm J}(0)>0$. By Lemma~\ref{lem:majvarodd}, it is enough for this to have ${\rm v}_{n,d} < (1-\upmu_{n,d}(0)) \uprho^2$. 
Setting $v={\rm v}_{n,d}$, and using the bounds ${\rm Var}\, I \leq v$ and $\upmu_{n,d}(0) \leq \sqrt{\frac{2}{\pi v}}$ of Lemmas~\ref{lem:majvarodd} and~\ref{lem:boundmu},  
the positivity of ${\rm E}_{n,d}-\upmu_{n,d}(0){\rm J}(0)$ is implied by
$$v < (1- \sqrt{\frac{2}{\pi v}})\uprho^2.$$
Observe that the map $f : \R_{>0} \rightarrow \R$, $x \mapsto (1- \sqrt{\frac{2}{\pi x}})\uprho^2 - x$, is 
continuous and strictly concave, with limits $-\infty$ at both $0$ and $+\infty$, and $f(8/\pi)=\uprho^2/2-8/\pi>0$. 
It has thus exactly two zeros $x_1 < x_2$, and a computation shows
$$ x_1 = 0.64484...\, \, {\rm and}\, \, x_2=92.43239...$$
See Figure~\ref{fig:jensen-even} below. Note that we have $v=\frac{n(d-2)}{3d}\geq 2/3>x_1$ for $n\geq 6$ or $d\geq 4$, so we have 
$f(v)>0$ if and only if $v<x_2$. We conclude assertion (ii) of the proposition by the numerical value 
$3\,x_2 \approx 277.29718...$
and we are done.\end{pf}

In the proof above we have used the following simple lemmas.
\begin{lemma} 
\label{lem:majvarodd} 
Assume $n,d\geq 3$. If $I$ is a $\Z$-valued random variable with law $\upmu_{n,d}$, we have 
${\rm Var}\, I \leq {\rm v}_{n,d}+2/3$ for $n$ odd, and ${\rm Var}\, I \leq {\rm v}_{n,d}$ for $n$ even.
\end{lemma}

\begin{pf} 
This is obvious by~\eqref{eq:recallvar} for $n$ even, so we may assume $n$ is odd. 
Setting $t=d-1$, we have $d\,\upbeta_{n,d}=(d-1)(t^{n-1}-1)$ by \eqref{eq:formbetand}.
Using $t^i \geq 2^i \geq i+1$ for all $i\geq 1$,
we deduce $\frac{t^{n-1}-1}{t-1} \geq n(n-1)/2$, hence $\frac{d}{d-1} \upbeta_{n,d}  \geq (d-2)n(n-1)/2$.
So we have ${\rm v}_{n,d}\,(n-1)(d-1)/ \upbeta_{n,d} \leq 2/3$ and we conclude by \eqref{eq:recallvar}.
\end{pf}

\begin{lemma} 
\label{lem:hnd2nonzero}
Let $n,d\geq 3$ be integers with $n$ even and $(n,d) \neq (4,3)$. Then we have $\upmu_{n,d}(2)>0$,
and in particular, $\upmu_{n,d}(0)<1$. 
\end{lemma}

\begin{pf} 
Fix $i\geq 0$ even. We have $\upmu_{n,d}(i)={\rm h}_{n,d}(i)/\upbeta_{n,d}$ 
with ${\rm h}_{n,d}(i)$ the number of ways to write the integer $d(n+i)/2$ 
as a sum of $n$ integers between $1$ and $d-1$. 
As $d(n+i)/2\geq 3n/2\geq n$, 
we have ${\rm h}_{n,d}(i)>0$ if and only if
$d(n+i)/2 \leq n(d-1)$. But this condition also reads $di \leq nd-2n$,
{\it i.e.} $2 i \leq (n-i)(d-2)$. 
So we have ${\rm h}_{n,d}(2)>0$ for $n \geq 6$ or $d \geq 4$.
\end{pf}

In the following lemma, we set ${\rm u}_{n,d}= 1+(-1)^n/(d-1)^{n-1}$.
Note that we have ${\rm u}_{n,d} \geq 1$ for $n$ even.

\begin{lemma} 
\label{lem:boundmu} 
For all $n,d \geq 3$, and $i \in \Z$, we have $\upmu_{n,d}(i) <  \frac{1}{{\rm u}_{n,d}}\sqrt{\frac{2}{\pi {\rm v}_{n,d}}}$. 
\end{lemma}

\begin{pf}
The main result of \cite{mattnerroos} (with $\ell = d-1 \geq 2$ and $n\geq 3$) 
asserts that for all $x \in \R$ we have 
$\mathds{P}(S_{n,d}=x) < {\rm M}_{n,d}$ with ${\rm M}_{n,d}=\sqrt{\frac{6}{\pi n d(d-2)}}$. 
It follows that for all $i \in \Z$ we have 
$\upmu_{n,d}(i) < {\rm M}_{n,d}/\mathds{P}(\mathcal{C}_{n,d})$. 
We conclude as 
$\mathds{P}(\mathcal{C}_{n,d})=(d-1)^{-n}\upbeta_{n,d} = {\rm u}_{n,d}/d$, 
by Formula~\eqref{eq:formbetand}.
\end{pf}

\begin{figure}[htbp]
\centering
\includegraphics[width=0.65\textwidth]{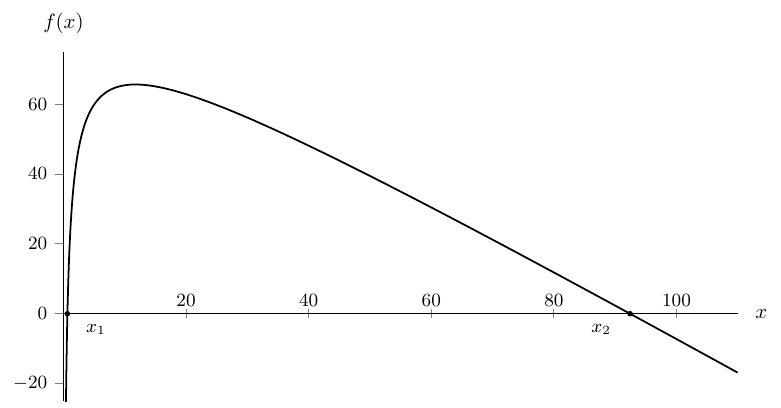}
\caption{The function
$f(x)=\uprho^2\bigl(1-\sqrt{\frac{2}{\pi x}}\bigr)-x$.}
\label{fig:jensen-even}
\end{figure}
  
\subsection{About the possible values of $n(d-2)/d$}

The two numerical constants appearing in Proposition~\ref{prop:inegjensen} 
are given with a slightly pedantic precision (especially since these are not the 
best bounds on $n(d-2)/d$ as we shall see). We did so because they do not weaken the corresponding exact numerical bounds, as explained by the following remark. 

\begin{remark}
\label{rem:possvaluesv}
{\it Assume $n,d\geq 3$ are integers. 
We have $n/3 \leq n(d-2)/d \leq n$. 
Moreover, for any real $C>0$ we have 
$$n(d-2)/d \leq C \, \, \iff\, \,{\rm either}\,\, n\leq C,\, \, {\rm or}\, \, n \in ]C,3C] \,\,{\rm and}\, \, d \leq 2n/(n-C).$$
In particular, there are only finitely many $(n,d)$ with $n>C$ and $n(d-2)/d \leq C$.  
}
\end{remark}

It follows that in any interval of $\R$ of the form $]a,b[$ with $a \in \N$ and $(b-a)<1$, there are only finitely many 
elements of the form $n(d-2)/d$, and a trivial algorithm to enumerate them all.\ps

\begin{example}
\label{ex:numgiveval} 
{\it We easily check that the last $n(d-2)/d$ before $3\uprho^2-2$ is $v_1=300.39253...$ for 
$(n,d)=(301,991)$, and that the first after $3\uprho^2-2$ is $v_2=300.39314...$ for $(n,d)=(301,992)$, 
explaining the optimality of the chosen constant $300.3928$ 
(but of course, any constant in $]v_1,v_2[$ works equally well).}
\end{example}

Such arguments will again play a role in Section~\ref{sect:pfmainthm} as well.

\subsection{Some conjectural generalizations of the Jensen bound}${}^{}$
\label{subsect:generaljensen}
The arguments in~\S\ref{subsect:mainexplform} and \S\ref{subsect:jensen} are sufficiently simple and general 
to apply, under standard conjectures, to general pure motives over $\Q$ rather than the middle degree cohomology of a smooth hypersurface. Our aim in this paragraph is to state the general predictions that they offer, {\it assuming the standard conjectures on the realizations and ${\rm L}$-functions of pure motives} (see e.g.~\cite{deligneperiodes, serre_facteurs_locaux}). \ps

Let $M$ be a nonzero pure motive over $\Q$ of weight $w \in \Z$ and with coefficients in $\Q$. 
We assume $w$ odd to simplify.
Define the {\it Hodge absolute first moment} of $M$ and the {\it Hodge variance} of $M$ respectively as 
\begin{equation}
\label{eq:ABSMVARM}
{\rm Abs}\,M \,=\,\sum_{p,q} |p-q| \frac{h^{p,q}(M)} {\dim M}\, \,\,\, {\rm and}\, \,\,\, {\rm Var}\, M \,=\, \sum_{p,q} (p-q)^2 \frac{h^{p,q}(M)}{\dim M},
\end{equation}
where $h^{p,q}(M)$ denote the Hodge numbers of (the de Rham realization of) $M$. 
The Cauchy-Schwarz inequality shows ${\rm Abs}\, M \,\leq\, \sqrt{ {\rm Var}\,M}$.
We denote by ${\rm Art}\,M$ the Artin conductor of (the $\ell$-adic realizations of) $M$, which is an integer $\geq 1$.
We finally denote by $\Lambda(s,M)$ the normalized complete ${\rm L}$-function of $M$.
Conjecturally, this is an entire function of $s \in \C$, bounded in all vertical strips, and satisfying $\Lambda(s,M)= \epsilon N^{1/2-s} \Lambda(1-s,M^\vee)$ for some $\epsilon \in \C^\star$ and $N={\rm Art}\, M$. 
The key inequality given by the Explicit Formula, together with the Jensen inequality and the fact that ${\rm J}$ is decreasing, simply reads: 
\begin{equation}
\label{eq:EFineqgen}
{\rm J}(\sqrt{{\rm Var}\,\,M}) \leq {\rm J}({\rm Abs}\,M) \leq  \frac{\log {\rm Art}\, M}{\dim M}.
\end{equation}

In particular, if $M$ has good reduction everywhere, {\it i.e. } ${\rm Art}\,M\,=1$, then we have ${\rm Var}\,M \geq \uprho^2$.
For instance, if we have furthermore $|p-q| \leq 10$ whenever $h^{p,q}(M) \neq 0$ then we have ${\rm Var}\,M \leq 10^2=100 < \uprho^2$ : so a nonzero $M$ should not exist. This corroborates for instance the Fontaine-Abrashkin theorem (case $M={\rm h}^1(A)$ with $A$ an abelian variety over $\Z$, with ${\rm Var}\, M=1$), or the fact that there is no cusp form of weight $<12$ for ${\rm SL}_2(\Z)$. 
This is also the main argument used in the proof of Proposition~\ref{prop:inegjensen} (i). Such arguments offer other conjectural perspectives. For instance, we may apply them to show the conjectural nonexistence of certain complete intersections in ${\mathbb P}^n_\C$ which have a proper smooth model over $\Z$ (using the Libgober-Wood variance formula \cite[Prop. 2.3]{LW} for their variance). We omit the details.

\section{The weighted Gaussian comparison}
\label{sect:weightedcomp}

\subsection{The main statement} 
\label{subsect:statementgaussiancomp}

As already explained in \S\ref{subsect:ideasandorg}, 
our basic strategy to prove Theorem~\ref{main:thm} 
is to compare the probability measure $\upmu_{n,d}$ with 
a (centered) Gaussian measure on $\Z$ with parameter ${\rm v}_{n,d}=\frac{n(d-2)}{3d}$.
We first introduce this measure. 
For $\epsilon \in \{0,1\}$ and $v>0$, we set 
$\Z_\epsilon = \epsilon + 2\Z$ and ${\rm m}_\epsilon(v)=\sum_{i \in \Z_\epsilon} e^{-i^2/2v}$.

\begin{definition} 
\label{def:gammaveps}
{\it 
Fix $\epsilon \in \{0,1\}$ and $v>0$. The {\rm Gaussian probability measure on $\Z_\epsilon$  
of parameter $v$} is the discrete probability $\gamma_v^{\epsilon}$ on $\Z$ with support $\Z_\epsilon$ satisfying 
$$\gamma_v^\epsilon(i) = \frac{e^{-i^2/2v}}{{\rm m}_\epsilon(v)} \, \, \, {\rm for\, \, all}\, \, \, i \in \Z_\epsilon.$$
}
\end{definition}

We shall not need it, but one could check using Lemma~\ref{lem:cargaussian} below that the variance of $\gamma_v^{\epsilon}$ is exponentially close to $v$ as $v \rightarrow \infty$. 
Our analysis will not only apply to the Odlyzko functions, but more generally to 
any bounded integrable function $H : \R_{\geq 0} \rightarrow \R$ that is differentiable at $0$. 
Such an $H$ will be called an {\it admissible weight}.  

\begin{definition}
\label{def:thetaepsvH}
{\it Let $\epsilon \in \{0,1\}$, $v>0$, and let $H$ be an admissible weight.
If $I$ denotes any discrete $\Z$-valued random variable with law $\gamma_v^\epsilon$, we set\footnote{
This sum is absolutely convergent, as ${\rm J}_H(|i|) = {\rm O}({\rm log}\,|i|\,)$ for $i \rightarrow +\infty$ 
by Lemma~\ref{lem:logcrois}.}
\begin{equation}
\label{eq:defthetaepsilonH}
\Theta_\epsilon^H(v):= \mathds{E}( \,{\rm J}_H(|I|)\,)\,=\,\sum_{i \in \Z} \,\gamma_v^\epsilon(i)\,{\rm J}_H(|i|).
\end{equation}}
\end{definition}

We now state the precise comparison we shall prove between ${\rm E}_{n,d}^H$ and 
$\Theta_\epsilon^H({\rm v}_{n,d})$, 
with explicit bounds that are uniform in $n$ and $d$. These bounds involve the following specific constants that will naturally enter the proof:
\begin{equation}\label{eq:twoconstants}  {\rm u}_{n,d} = 1+(-1)^n/(d-1)^{n-1}\, \, \, {\rm and}\, \, \, \eta_d = \frac{(d-1)^4 -1}{180 d^4}.
\end{equation}
We finally define the constant ${\rm C}_\epsilon$ by ${\rm C}_0=\pi/4$ and ${\rm C}_1=1/2$.

\begin{thm} 
\label{thm:weightedgaussian}
Assume $n\geq 5$, $d\geq 3$, $\epsilon \equiv n \bmod 2$ and set $v={\rm v}_{n,d}$. 
\begin{itemize}
\item[(i)] For any admissible weight $H$, we have the inequality
$$| {\rm E}_{n,d}^H - \Theta_\epsilon^H(v)|\, \leq \,\,{\rm C}_\epsilon ||H||_\infty \, ({\rm M}_{n,d}+{\rm R}_{n,d}),$$ 
with  ${\rm M}_{n,d}  = \frac{n}{{\rm u}_{n,d}}(\frac{4\,\eta_d}{v^2}+\frac{64}{8505\, v^3})$ and ${\rm R}_{n,d} \,=\, 17 n (7/8)^n$.\ps
\item[(ii)] If furthermore $H\geq 0$, we also have ${\rm E}_{n,d}^H - \Theta_\epsilon^H(v) \,\leq \,{\rm C}_\epsilon\, ||H||_\infty \, {\rm R}_{n,d}$.
\end{itemize}
\end{thm}

We refer to Section~\ref{sect:pfmainthm} for how to deduce Theorem~\ref{main:thm} from this result. 
Our goal in this section is to prove Theorem~\ref{thm:weightedgaussian}.

\subsection{The Fourier reduction}${}^{}$
\label{subsect:fourierreduction}
	In this section, we fix $\epsilon \in \{0,1\}$ and any signed, discrete, measure 
$\mu$ on $\Z$ 
with total mass $\mu(\Z)=0$, support in $\Z_\epsilon$, and finite first moment $\sum_{j \in \Z} |j \mu(j)| <\infty$. 
Our motivation is the example $\mu=\upmu_{n,d}-\gamma^\epsilon_v$ with $v={\rm v}_{n,d}$ and $\epsilon \equiv n \bmod 2$, which satisfies
${\rm E}_{n,d}^H - \Theta_\epsilon^H(v) = \sum_{j \in \Z_\epsilon} \mu(j) {\rm J}_H(|j|)$, 
a natural starting point for proving Theorem~\ref{thm:weightedgaussian}.\ps

\begin{lemma} 
\label{lem:signedmeasure}
For any admissible weight $H$, we have the equality 
$$\sum_{j \in \Z} \mu(j) {\rm J}_H(|j|) \,= \, \int_0^\infty\frac{H(t)}{2\,{\rm sinh}(t/2)}  
\left(\sum_{j \in \Z_\epsilon} \mu(j) e^{-|j|t/2}\right) dt,$$
where both the sum and the integral above are absolutely convergent.
\end{lemma}

\begin{pf} 
The left-hand side is absolutely convergent as we have 
$|{\rm J}_H(|j|)| = {\rm O}({\rm log}\,|j|)$ for $|j| \rightarrow \infty$. 
We have $\sum_{j \in \Z} \mu(j)=0$ (absolutely convergent sum) by assumption,
so the left-hand side coincides with $\sum_{j \in \Z} \,\mu(j)\, ({\rm J}_H(|j|)-{\rm J}_H(0))$.
But by definition of ${\rm J}_H$ we have, for all $j \in \Z$, 
$${\rm J}_H(|j|)-{\rm J}_H(0)\,=\,\int_0^\infty H(t)e^{-t/2}(e^{-|j|t /2}-1)/(1-e^{-t})dt.$$
Setting $h_{j}(t)=H(t)e^{-t/2}(e^{-|j|t/2}-1)/(1-e^{-t})$, and using $1-e^{-x}\leq x$ for $x \geq 0$,
we have $|h_{j}(t)| \leq ||H||_\infty |j| \frac{t}{4 {\rm sinh}\,(t/2)}$, 
hence $\sum_{j \in \Z} \int_0^\infty |\mu(j)||h_j(t)|dt < \infty$ by 
the assumption on $\mu$ and $\int_0^\infty t/({\rm sinh} t) dt <\infty$. 
But for all $t>0$ we have 
$$\sum_{j \in \Z} \mu(j) h_{j}(t) \,=\,\frac{H(t)}{2 \, {\rm sinh\,} (t/2)}\, \sum_{j \in \Z} \mu(j) e^{-|j|t/2},$$
and we conclude by Fubini's theorem.
\end{pf}

Our aim now, standard in approaches to the local central limit theorem, is to apply a Fourier inversion on the circle to express the sum in parentheses above. If $\lambda$ is any discrete, signed measure on $\Z$ with $\sum_j |\lambda(j)|<\infty$, and supported in $\Z_\epsilon$, recall that its {\it Fourier transform} (or {\it characteristic function}) is the absolutely convergent Fourier series
$$\widehat{\lambda}(x) = \sum_{j \in \Z} \lambda(j) e^{i j x}, \, \, \, x \in \R.$$
This is a continuous $2\pi$-periodic function on $\R$, 
satisfying $\widehat{\lambda}(x+\pi) = (-1)^\epsilon \widehat{\lambda}(x)$ for all $x \in \R$. 
Of course, $\lambda=\mu$ is an example, with $\widehat{\mu}(x)$ of class $\mathcal{C}^1$ and satisfying $\widehat{\mu}(0)=\mu(\Z)=0$.
Another example is given, for $t>0$, by $\lambda=\lambda_t^\epsilon$, with $\lambda_t^\epsilon(j) = e^{-|j| t/2}$ 
for all $j \in \Z_\epsilon$ and $0$ otherwise. It satisfies $\widehat{\lambda_t^\epsilon}(x) = \sum_{j \in \Z_\epsilon} e^{-|j|t/2} e^{i j x}$. 

\begin{lemma} 
\label{lem:0fourier}For all $t>0$ we have $\sum_{j \in \Z_\epsilon} \mu(j) e^{-|j| t/2} = \frac{1}{\pi} \int_{-\pi/2}^{\pi/2} \widehat{\mu}(x) \widehat{\lambda_t^\epsilon}(-x) dx$.
\end{lemma}

\begin{pf} The product of the two absolutely convergent Fourier series $\widehat{\mu}(x)$ and $ \widehat{\lambda^\epsilon_t}(-x)$ is $\pi$-periodic, with constant Fourier coefficient equal to the left-hand side above.
\end{pf}

The Fourier series $\widehat{\lambda_t^\epsilon}(x)$ has a simpler and useful expression.

\begin{lemma}
\label{lem:expressionlambdatepshat} Fix $t>0$. We have the formulas $\widehat{\lambda_t^0}(x) = \frac{{\rm sinh}(t)}{{\rm cosh}(t) - {\rm cos}(2x)}$ and $\widehat{\lambda_t^1}(x) = \frac{2\, {\rm sinh}(t/2) \,{\rm cos}(x)}{{\rm cosh}(t) - {\rm cos}(2x)}$. In particular, we have $\widehat{\lambda_t^\epsilon}(x)>0$ for $|x|<\pi/2$.
\end{lemma}

\begin{pf} 
Set $z=e^{i x} e^{-t/2}$, so that $|z|<1$. A geometric summation shows thus 
$$\widehat{\lambda_t^1}(x) \,=\, {\rm Re} \,\frac{2z}{1-z^2}\, \, \, {\rm and}\, \, \widehat{\lambda_t^0}(x) \,=\, 1\,+\, 2 \,{\rm Re}\, \frac{z^2}{1-z^2} \,=\, {\rm Re}\, \frac{1+z^2}{1-z^2}.$$
Note first that we have $|1-z^2|^2\, =\, 1+|z|^4-2 {\rm Re}\, z^2 \,=\, 2e^{-t}( {\rm cosh}(t) - {\rm cos}(2x) )$ and 
${\rm Re}\, 2z(1-\overline{z}^2) \,= \,2 (1-|z|^2) \,{\rm Re}\, z \,=\, 4 e^{-t} {\rm sinh}(t/2) {\rm cos}(x)$, hence the formula 
for $\widehat{\lambda_t^1}$. The formula for $\widehat{\lambda_t^0}$ follows similarly using the identity
${\rm Re}\, (1+z^2)(1-\overline{z}^2) \,=\, 1-|z|^4 + 0 \ \,=\, 2 e^{-t}\, {\rm sinh}(t)$.
\end{pf}
	
This suggests to introduce the following notation. 
Let us denote by $\mathcal{F}$ the space of functions $\phi : [-\pi/2,\pi/2] \rightarrow \C$ 
of class $\mathcal{C}^1$ satisfying $\phi(0)=0$. 
For an admissible weight $H$, we finally define 
a linear form ${\rm T}_\epsilon^H : \mathcal{F} \rightarrow \C$ by the formula 
\begin{equation}
\label{eq:defTepsH} 
{\rm T}_\epsilon^H(\phi)=\int_0^{\infty}\int_{-\pi/2}^{\pi/2}  \frac{\widehat{\lambda^\epsilon_t}(-x)}{2\pi\, {\rm sinh}(t/2)} \phi(x) H(t)\, {\rm d}x\, {\rm d}t.
\end{equation}
As we shall see in the proof of Prop.~\ref{prop:normTHeps} below, this integral is absolutely convergent for $\phi \in \mathcal{F}$. Lemmas~\ref{lem:signedmeasure} and~\ref{lem:0fourier} imply:

\begin{cor} 
\label{cor:newfrmula}
For any admissible weight $H$, we have $\sum_{j \in \Z} \mu(j) {\rm J}_H(|j|) \,=\, {\rm T}_\epsilon^H(\widehat{\mu})$. 
\end{cor}

The formula  $||f||_\star=\int_{-\pi/2}^{\pi/2} \,|f(x)/x|\, {\rm d}x$
defines a norm on $\mathcal{F}$. Recall the constants ${\rm C}_0=\pi/4$ and ${\rm C}_1=1/2$ defined in \S\ref{subsect:statementgaussiancomp}.
A key estimate is the following.

\begin{prop} 
\label{prop:normTHeps}
For any $\phi \in \mathcal{F}$, and any admissible weight $H$, we have:
\begin{itemize}
\item[(i)]  $|{\rm T}_\epsilon^H(\phi)| \leq {\rm C}_\epsilon ||H||_\infty ||\phi||_\star$,\ps
\item[(ii)] ${\rm T}_\epsilon^H(\phi) \geq 0$ for $\phi\geq 0$ and $H \geq 0$.
\end{itemize}
\end{prop}

\begin{pf} 
Recall $\widehat{\lambda_t^\epsilon}(x) \geq 0$ for $|x| \leq \pi/2$ 
by Lemma~\ref{lem:expressionlambdatepshat}.
Set ${\rm K}^\epsilon(x)= \int_0^\infty \frac{\widehat{\lambda_t^\epsilon}(x)}{ 2 \pi \,{\rm sinh}(t/2)} dt$.
We claim ${\rm K}^\epsilon(x) \leq {\rm C}_\epsilon/|x|$ for $0<|x| < \pi/2$. 
This claim immediately implies
$$\int_{-\pi/2}^{\pi/2}  \int_0^\infty \frac{\widehat{\lambda^\epsilon_t}(-x)}{2\pi\, {\rm sinh}(t/2)}\, |\phi(x)|\,|H(t)|\, {\rm d}t \,{\rm d}x \leq {\rm C}_\epsilon ||H||_\infty ||\phi||_\star,$$
hence both the absolute convergence of \eqref{eq:defTepsH} by Fubini, and (i) and (ii). \par
We now prove the claim for $\epsilon=0$. 
By Lemma~\ref{lem:expressionlambdatepshat}, we have 
$$\widehat{\lambda_t^0}(x)=\frac{{\rm sinh}(t)}{{\rm cosh}(t)-{\rm cos}(2x)}=\frac{2\,{\rm sinh}(t/2)\,{\rm cosh}(t/2)}{2\,{\rm sinh}(t/2)^2 + 2 \,{\rm sin}(x)^2}.$$
For $x \notin \pi \Z$, the change of variable $u=\frac{{\rm sinh}(t/2)}{|{\rm sin}\, x|}$ 
shows ${\rm K}^0(x)=\frac{1}{\pi |{\rm sin} x|} \int_0^\infty \frac{1}{u^2+1} du = \frac{1}{2 | {\rm sin}\,x|}$. 
We conclude by the convexity inequality $|{\rm sin}\, x| \geq 2|x|/\pi$ for $|x| \leq \pi/2$.\ps
We finally assume $\epsilon =1$. Lemma~\ref{lem:expressionlambdatepshat} shows 
$\widehat{\lambda^1_t}(x) = \widehat{\lambda^0_t}(x) \frac{{\rm cos}\,x} { {\rm cosh}(t/2)} \leq {\rm cos}(x) \widehat{\lambda^0_t}(x)$, 
so ${\rm K}^1(x) \,\leq \,{\rm K}^0(x) \,{\rm cos}(x) = \frac{1}{2 |{\rm tan}\, x|}$ by the above computation. 
We conclude by the convexity inequality $|{\rm tan}\,x| \geq |x|$ for $|x| < \pi/2$.
\end{pf}

\subsection{The two characteristic functions}
\label{subsect:charfun}

In this paragraph, we determine the two characteristic functions $\widehat{\upmu_{n,d}}$ and $\widehat{\gamma_v^\epsilon}$. Recall the quantity ${\rm u}_{n,d} = 1 + (-1)^n /(d-1)^{n-1}$ and consider, for $x \in \R \smallsetminus d\pi \Z$,
 $$\varphi_d(x) = \frac{\sin\, (d-1)x/d}{(d-1)\sin x/d},$$
a function that extends continuously to $\varphi_d : \R \rightarrow \R$.
\begin{lemma} 
\label{lem:carhodge}
For all integers $n,d \geq 3$ and all $x \in \R$ we have $$\widehat{\upmu_{n,d}}(x) = \frac{1}{{\rm u}_{n,d}}\sum_{k=0}^{d-1} (-1)^{nk} \varphi_d(x+k \pi)^n.$$
\end{lemma}

\begin{pf} 
Recall from the proof of Lemma~\ref{lem:griffitsproba} that 
the generating series of the random variable $d S_{n,d}$ is ${\rm P}_d(u)^n$, 
with ${\rm P}_d(u)$ the Laurent polynomial defined in~\eqref{eq:defpolPd}. 
In other words, we have 
$\mathds{E}( e^{ix\, d S_{n,d}}) = {\rm P}_d(e^{ix})^n = \varphi_d(dx)^n$, 
{\it i.e.} $\mathds{E}( e^{ix\,S_{n,d}}) =\varphi_d(x)^n$, for all $x \in \R$.
The law $\upmu_{n,d}$ is that of $S_{n,d}$ conditioned to the event 
$S_{n,d} \in n+2\Z$, called $\mathcal{C}_{n,d}$ {\it loc. cit.} 
Note that $d S_{n,d}$ trivially takes its values in $nd + 2 \Z$, 
as well as the equality $ \mathds{P}(\mathcal{C}_{n,d}) = {\rm u}_{n,d}/d$ 
proved in Lemma~\ref{lem:griffitsproba}. 
The statement follows then from the fact that for $m \in nd+2\Z$ we have
$\frac{1}{d} \sum_{k=0}^{d-1} (-1)^{nk} e^{i m k\pi / d} = 1$ 
if $m \in d(n +2\Z)$ and $0$ otherwise.
\end{pf}

\begin{remark} 
\label{rem:s2d}
{\it Assume $d\geq 3$. The random variable $S_{2,d}$ of \S\ref{subsect:Snd} still makes sense. It satisfies $S_{2,d} \in ]-2,2[$, and thus $S_{2,d} \in 2\Z \Leftrightarrow S_{2,d}=0$, and then $\mathds{P}(\mathcal{C}_{2,d})=1/(d-1)>0$. We can thus still define $\upmu_{2,d}$ as the law of $S_{2,d}$ given $\mathcal{C}_{2,d}$, and it satisfies $\widehat{\upmu_{2,d}}(x)=1$ identically. The proof of Lemma~\ref{lem:carhodge} applies verbatim in this case, with ${\rm u}_{2,d} \,=\, d \,\mathds{P}(\mathcal{C}_{2,d})= d/(d-1)$, and shows that for all $x \in \R$ we have $\frac{d}{d-1} \, = \,\sum_{k=0}^{d-1} \varphi_d(x+k\pi)^2$. We will use this identity below.}
\end{remark}

Fix now $\epsilon \in \{0,1\}$ and $v>0$, and recall the mass ${\rm m}_\epsilon(v)=\sum_{j \in \Z_\epsilon} e^{-j^2/2v}$.
The characteristic function of the Gaussian measure $\gamma_v^\epsilon$ is by definition
$$\widehat{\gamma_v^\epsilon}(x) \, =\,\frac{1}{{\rm m}_\epsilon(v)} \,\sum_{j \in \Z_\epsilon} \,e^{-j^2/{2v}}\, e^{i j x}.$$

\begin{lemma} 
\label{lem:cargaussian}
For all $x \in \R$ we have $\widehat{\gamma_v^\epsilon}(x) \,=\, \frac{1}{{\rm M}_\epsilon(v)} \sum_{k \in \Z} (-1)^{k \epsilon} e^{-v( x - \pi k )^2/2}$, with ${\rm M}_\epsilon(v)\,=\,2\,{\rm m}_\epsilon(v)/\sqrt{2 \,\pi v}\, = \sum_{k \in \Z} (-1)^{k \epsilon} e^{-v\pi^2 k^2/2}$.
\end{lemma}

\begin{pf} 
The function $f_{v,\epsilon} : \R \rightarrow \C, u \mapsto e^{-(2u+\epsilon)^2/{2v}} e^{i(2u+\epsilon)x}$ 
is a Schwartz function on $\R$. The standard Fourier transform of the Gaussian, 
as well as the change of variables $s=2u+\epsilon$, show that for all $\xi \in \R$ we have
{\small
$$\int_\R f_{v,\epsilon}(u) e^{-2i\pi u \xi} du \,=\, \frac{1}{2}\, \int_\R \,e^{-s^2/2v} \,e^{i s x -  i \pi \xi (s- \epsilon) }\,ds\,=\, \frac{1}{2}\,e^{i \pi \xi \epsilon}\,\,\sqrt{2 \pi v} \,\,e^{-v(x-\pi \xi)^2/2}.$$
}
\par \noindent The Poisson formula for $f_{v,\epsilon}$ is the formula of the statement, 
with ${\rm M}_\epsilon(v)\,=\,2\,{\rm m}_\epsilon(v)/\sqrt{2 \pi v}>0$. 
The last formula for ${\rm M}_\epsilon(v)$ follows from $\widehat{\gamma_v^\epsilon}(0)=1$.
\end{pf}

\subsection{The main term}

Recall the function $\varphi_d : \R \rightarrow \R$ introduced in \S\ref{subsect:charfun}.
We now compare $\varphi_d(x)^n$ with $e^{-vx^2/2}$, where $v={\rm v}_{n,d}$. 
The resulting bound for the norm $||\cdot||_\star$ of their difference will give rise to the main term ${\rm M}_{n,d}$ in Theorem~\ref{thm:weightedgaussian}. As a preliminary step, we begin with the following expansion of $\varphi_d$, which explains the appearance of the constant $\eta_d$.

\begin{lemma} 
\label{lem:boundcentral}
For $d\geq 2$ and $|x| \leq \pi/2$ we have 
$$\varphi_d(x)  = e^{ - \frac{d-2}{3d} \frac{x^2}{2} - \eta_d x^4 - {\rm R}_d(x)}, \, \, \, {\rm with}\,\, \,0 \leq {\rm R}_d(x) \leq \frac{4}{8505} x^6.$$
\end{lemma}

\begin{pf} 
For $0<|x|<\pi$, we have $(\sin \,x)/x> 0$, so by Euler's product formula,
\begin{equation}
\label{eq:EulerE}
{\rm E}(x)  \,:=\, - \log \frac{\sin x}{x}\, = \, \sum_{m \geq 1} - \log (1- \frac{x^2}{ \pi^2 m^2}) = \sum_{k \geq 1} \frac{\zeta(2k)}{k \,\pi^{2k}} x^{2k}.
\end{equation}
We also have $\varphi_d(x) = \frac{ {\rm sin}\, (d-1)x/d}{(d-1)x/d}\frac{x/d}{{\rm sin}\, x/d}$, and both fractions are $>0$, so 
$$- \log \varphi_d(x) \,=\, {\rm E}( (d-1)x/d ) - {\rm E}(x/d)\, =\, \sum_{k\geq 1}  \,a_{2k} \,x^{2k},$$
with $a_{2k} = \frac{\zeta(2k)}{\pi^{2k}\,k} ((\frac{d-1}{d})^{2k}-\frac{1}{d^{2k}})$. 
Using $\zeta(2) = \pi^2/6$ and $\zeta(4)= \pi^4/90$ we find $a_2 = \frac{1}{6} (d-2)/d$ and $a_4=\eta_d$.
Using $0 \leq a_{2k} \leq \frac{1}{3} \frac{\zeta(6)}{\pi^{2k}}$ for $k\geq 3$, and setting ${\rm R}_d(x)= \sum_{k\geq 3}  \,a_{2k} \,x^{2k}$, we obtain 
$0 \leq {\rm R}_d(x) \leq \frac{x^6}{3} \frac{\zeta(6)}{\pi^6} \frac{1}{1-x^2/\pi^2}$.
Using finally $\zeta(6)= \pi^6/945$, we obtain for $|x| \leq \pi/2$ the inequality 
$0 \leq {\rm R}_d(x) \leq \frac{4}{8505} x^6$.
\end{pf}

\begin{lemma} 
\label{lem:boundcentralterm}
For $n,d \geq 3$ and $v={\rm v}_{n,d}$, we have : \begin{itemize}
\item[(i)] $\varphi_d(x)^n - e^{- v x^2/2} \leq 0$ for $|x| \leq \pi/2$, \ps
\item[(ii)] $||\varphi_d(x)^n-e^{- v x^2/2}||_\star \leq 4\eta_d \frac{n}{v^2} +\frac{64}{8505} \frac{n}{v^3}$.
\end{itemize} 
\end{lemma}

\begin{pf} 
By Lemma~\ref{lem:boundcentral} and $0 \leq 1-e^{-u} \leq u$ for $u\geq 0$, we have for $|x| \leq \pi/2$
$$0 \leq e^{- v x^2/2} - \varphi_d(x)^n = e^{ - v \frac{x^2}{2}}(1- e^{- n \eta_d x^4 - n {\rm R}_d(x) })  \leq n e^{ - v \frac{x^2}{2}} (\eta_d x^4 + \frac{4}{8505} x^6).$$
This shows (i). In order to show (ii), we apply $||.||_\star$ and deduce 
$$||\varphi_d(x)^n-e^{- v x^2/2}||_\star \,\leq \,n \,(\eta_d \,{\rm I}_4\,+ \,\frac{4}{8505} \,{\rm I}_6)$$
with ${\rm I}_k =|| e^{-vx^2/2} x^k ||_\star$ for $k\geq 1$. We conclude by 
the inequality ${\rm I}_k \leq \int_\R |x|^{k} e^{-v x^2/2} \frac{dx}{|x|}$, this last integral being  
$\Gamma(k/2)(2/v)^{k/2}$ by Lemma~\ref{lem:calcgammaint}.
\end{pf}

During the proof, we have used the following formula.

\begin{lemma} 
\label{lem:calcgammaint}
For $u>0$ and $m>0$ we have 
$\int_\R |x|^{2m} e^{-u x^2} \frac{dx}{|x|}= \frac{\Gamma(m)}{u^m}$.
\end{lemma}

\begin{pf} 
Set $t = ux^2$ in $\Gamma(m) = \int_0^\infty t^m e^{-t} dt/t$.
\end{pf}

\subsection{The remaining terms}

For $n,d \geq 3$ we set $\Phi_{n,d}(x) = \sum_{k=1}^{d-1} (-1)^{kn} \varphi_d(x+k \pi)^n$.
By Lemma~\ref{lem:carhodge} we have 
$$ \widehat{\upmu_{n,d}}(x) = \frac{1}{{\rm u}_{n,d}} ( \varphi_d(x)^n + \Phi_{n,d}(x)).$$
Similarly, for $v>0$, $\epsilon \in \{0,1\}$ and ${\rm G}^\star_v(x)=\sum_{k \neq 0} (-1)^{k\epsilon} e^{- v (x-k\pi)^2/2}$, we have
$$\widehat{\gamma_v^{\epsilon}}(x) = \frac{1}{{\rm M}_\epsilon(v)}( e^{ - vx^2/2} + {\rm G}^\star_v(x))$$
by Lemma~\ref{lem:cargaussian}. Using $\widehat{\upmu_{n,d}}(0)=\widehat{\gamma_v^{\epsilon}}(0)=\varphi_d(0)=1$ we may thus write
\begin{multline}
\label{eq:decompfourier}
\widehat{\upmu_{n,d}}(x)-\widehat{\gamma_v^\epsilon}(x) = \frac{1}{{\rm u}_{n,d}} \bigl(\varphi_d(x)^n-e^{-vx^2/2}\bigr)  + \frac{1}{{\rm u}_{n,d}} \bigl(\Phi_{n,d}(x)-\Phi_{n,d}(0)\bigr) \\
- \frac{1}{{\rm M}_\epsilon(v)} \bigl({\rm G}^\star_v(x)-{\rm G}^\star_v(0)\bigr) + \left(\frac{1}{{\rm u}_{n,d}}-\frac{1}{{\rm M}_\epsilon(v)}\right) \bigl(e^{-vx^2/2}-1\bigr).
  \end{multline}
The four terms on the right belong to $\mathcal{F}$ and we want to bound their $||.||_\star$.
For $v={\rm v}_{n,d}$, the first term has been dealt with in Lemma~\ref{lem:boundcentralterm}. 
For the three remaining ones, we will repeatedly use the following trivial lemma.

\begin{lemma} 
\label{lem:majnormstar}
For $\phi \in \mathcal{F}$ and $M= {\rm sup}_{|x|\leq \pi/2} |\phi'(x)|$ we have $||\phi||_\star \leq \pi M$.
\end{lemma}

\begin{pf} 
We have $|\phi(x)| \leq |x| M$ by the mean value theorem and $\phi(0)=0$.
\end{pf} 

\begin{lemma} 
\label{lem:boundPhi}
For $n,d \geq 3$ we have 
$||\Phi_{n,d}(x)-\Phi_{n,d}(0)||_\star \leq 2 \pi n (\frac{\sqrt{3}}{2})^{n-1}$.
\end{lemma}

\begin{pf}
For all $x \in \R$ we have $(\varphi_d^n)'(x) = n \varphi_d'(x) \varphi_d(x)^{n-1}$, as well as
$|\varphi_d'(x)| \leq \mathds{E}(\frac{|X_1^{(d)}|}{d} )\leq 1$. By Lemma~\ref{lem:majnormstar}, it is thus enough to show
\begin{equation}
\label{eq:boundsumdec}
\sum_{k=1}^{d-1} |\varphi_d(x+k \pi)|^{n-1} \leq  2 (\frac{\sqrt{3}}{2})^{n-1}, \, \, \,{\rm for}\, \, |x| \leq \pi/2.
\end{equation}
On the one hand, Remark~\ref{rem:s2d} shows $\sum_{k=1}^{d-1} |\varphi_d(x+k \pi)|^{2} \leq  \frac{d}{d-1} \leq 3/2$. On the other hand, for $|x|\leq \pi/2$ and $1 \leq k \leq d-1$ we have 
$\pi/2d \leq (x+k\pi)/d \leq \pi - \pi/2d$, so Lemma~\ref{lem:boundsin} applied to $m=d-1$ shows $|\varphi_d(x+k\pi) |\leq \sqrt{3}/2$.
The left-hand side in~\eqref{eq:boundsumdec} is thus bounded above 
by $3/2 \cdot (\sqrt{3}/2)^{n-3} = 2 (\sqrt{3}/2)^{n-1}$. 
\end{pf}

We have used the following elementary bound.

\begin{lemma} 
\label{lem:boundsin}
For any integer $m\geq 2$ and $x \in [\frac{\pi}{2(m+1)},\pi-\frac{\pi}{2(m+1)}]$ we have $|\frac{\sin (mx)}{m \sin x} |\leq \sqrt{3}/2$. 
\end{lemma}

\begin{pf} The quantity $|\frac{\sin (mx)}{m \sin x}|$ and the interval of the statement are invariant under $x \mapsto \pi -x$, so we may assume $|x| \leq \pi/2$. Moreover, for $m=2$ we have $\frac{\sin (mx)}{m \sin x}={\rm cos} \,x$, which is maximal for 
$x=\frac{\pi}{2(m+1)}=\pi/6$ and is $\sqrt{3}/2$, so we may assume $m\geq 3$. Finally, for $\frac{\pi}{m} \leq x \leq \frac{\pi}{2}$, we have $|\sin x| \geq \sin \frac{\pi}{m} \geq \frac{2}{\pi} \frac{\pi}{m} =\frac{2}{m}$ (by concavity of $\sin$), and thus $|\frac{\sin (mx)}{m \sin x} |\leq 1/2<\sqrt{3}/2$. So we may assume 
$\frac{\pi}{2(m+1)} \leq x < \frac{\pi}{m}$.  But more generally for $0<x< \pi/m$, we have the bound
$$ 0 \leq \frac{{\rm sin}\, (mx)}{m \sin x} = e^{- ( {\rm E}(mx) - {\rm E}(x)) } \leq e^{ -(m^2-1)x^2/6},$$
which follows from Euler's formula~\eqref{eq:EulerE}. But for $x= \pi/(2(m+1))$ the right-hand side is 
$e^{ -\frac{m-1}{m+1}\pi^2/24}$, and this constant decreases with $m$ and is $e^{ -\pi^2/48} = 0.814... < \sqrt{3}/2$ for $m=3$. 
\end{pf}

\begin{lemma} 
\label{lem:boundGstar}
For $v\geq 4/\pi^2$ we have $||{\rm G}^\star_v(x)-{\rm G}^\star_v(0)||_\star \leq \, 2\, \pi^2 v \,e^{-\pi^2 v/8}.$
\end{lemma}

\begin{pf} 
By Lemma~\ref{lem:majnormstar}, 
it is enough to show $|({\rm G}^\star_v)'(x)| \leq 2 \pi v e^{-\pi^2 v/8}$ for all $|x| \leq \pi/2$. 
Consider the function $f : \R \rightarrow \R$, $y \mapsto y e^{-s y^2/2},$ with $s=v\pi^2$. 
Note that we have $|x/\pi-k|\geq |k|-1/2 \geq 1/2$ for $|x|\leq \pi/2$ and $k \neq 0$. 
Moreover, we have $f'(y) = (1-sy^2) e^{-sy^2/2}$ so $f$ is decreasing on $[1/2,+\infty[$
by the assumption $s\geq 4$. It follows that we have the inequalities
{\small
$$|({\rm G}^\star_v)'(x)| \,\leq \, \sum_{k\neq 0} \pi v f(|x/\pi-k|) \leq 2 \pi v \sum_{k=1}^{\infty} f(k-1/2) \leq 2\pi v \left( f(1/2) + \int_{\frac{1}{2}}^\infty f(y) dy\right).$$
}
We have $f(1/2)=e^{-s/8}/2$, and the integral on the right is exactly $e^{-s/8}/s \leq e^{-s/8}/4$, 
so the right-hand side above is $\leq 2 \pi v e^{-\pi^2 v/8}$.
\end{pf}

\begin{lemma} 
\label{lem:lastterm}
For $v>0$ we have $|| e^{-v x^2/2} - 1||_\star \leq v \pi^2/2$.
\end{lemma}

\begin{pf} 
Set $\phi(x)=e^{-v x^2/2} - 1$. 
For $|x| \leq \pi/2$ we have $|\phi'(x)|=v |x| e^{-v x^2/2} \leq v \pi/2$, 
so we conclude by Lemma~\ref{lem:majnormstar}.
\end{pf}

Our last task is to estimate the constants ${\rm u}_{n,d}$ and ${\rm M}_\epsilon(v)$ (that are exponentially close to $1$ when $n$ or $v$ grows). We gather these estimates in a final lemma.

\begin{lemma}
\label{lem:boundcst} For $n\geq 5$, $d\geq 3$ and $v={\rm v}_{n,d}$, we have:
\begin{itemize}
\item[(i)] ${\rm u}_{n,d} \geq 15/16$,\ps
\item[(ii)] ${\rm M}_\epsilon(v) \geq 4/5$,\ps
\item[(iii)] $|{\rm u}_{n,d}^{-1}-{\rm M}_\epsilon(v)^{-1}|\,\leq \,3 \,e^{-n \pi^2/18}+ 2^{2-n}$.
\end{itemize}
\end{lemma}

\begin{pf} 
We have $|{\rm u}_{n,d}-1|=1/(d-1)^{n-1} \leq 1/2^{n-1} \leq 1/16$ for $d\geq 3$ and $n\geq 5$, 
proving assertion (i). For (ii), we have 
\begin{equation}
\label{eq:majMepsv}
|{\rm M}_\epsilon(v) -1| \leq 2 \sum_{k=1}^\infty e^{-v \pi^2 k^2/2}.
\end{equation}
Set $q=e^{-\pi^2/18}$. Using $v\geq n/9$ and $k^2 \geq k$ we deduce 
$|{\rm M}_\epsilon(v) -1| \leq 2 q^n/(1-q^n)$ by a geometric summation. 
This upper bound decreases with $n$, and for $n=5$ it is $0.137... < 1/5$, 
hence (ii). For (iii), the bounds established above show
$$|{\rm u}_{n,d}^{-1}-{\rm M}_\epsilon(v)^{-1}| =\frac{ |({\rm M}_\epsilon(v)-1)-({\rm u}_{n,d}-1)| }{{\rm u}_{n,d} {\rm M}_\epsilon(v)} \leq \frac{2 \,q^n/(1-q^5) + 1/2^{n-1}}{3/4}$$
for $n\geq 5$. We conclude as $8/(3(1-q^5))=2.85...<3$.
\end{pf}

All the upper bounds above can be gathered into a simpler, but equally useful, bound.

\begin{lemma} 
\label{lem:concattail}
For $n\geq 5$, $d\geq 3$, $\epsilon \equiv n \bmod 2$ and $v={\rm v}_{n,d}$, we have 
$$\frac{2 \pi n}{{\rm u}_{n,d}} \, (\frac{\sqrt{3}}{2})^{n-1} + \frac{2 \pi^2 v}{{\rm M}_\epsilon(v)}\, e^{-\pi^2v/8} \,+\, |{\rm u}_{n,d}^{-1}-{\rm M}_\epsilon(v)^{-1}| \,\frac{v \pi^2}{2} \,\,< \,\, 17 \,n\, (7/8)^n.$$
\end{lemma}
\begin{pf} Using the three bounds in Lemma~\ref{lem:boundcst}, as well as the obvious inequalities $n/9 \leq v\leq n/3$, we are reduced to proving, for all $n\geq 5$, the inequality $$\frac{2\pi}{\alpha} \frac{16}{15} \alpha^{n}\, +\, \frac{1}{3}\,\frac{5}{4} \,2 \pi^2 \, \beta^n + \frac{1}{3}\,\frac{\pi^2}{2} \,(3\, \gamma^n\,+\, 4 \,\delta^n ) \,< 17 \,(7/8)^n,$$
with $\alpha = \frac{\sqrt{3}}{2} = 0.866...$, $\beta=e^{-\pi^2/72}=0.871...$, $\gamma = e^{-\pi^2/18}=0.577...$ and $\delta = 1/2$. But for $\tau=7/8 = 0.875$, we have $0<\alpha,\beta,\gamma,\delta < \tau$, and we conclude by the numerical computation
$$\frac{2\pi}{\alpha} \,\frac{16}{15} \,(\frac{\alpha}{\tau})^5 \, +\, \frac{1}{3}\, \frac{5}{4} \,2 \pi^2 (\frac{\beta}{\tau})^5\,+ \frac{1}{3}\,\frac{\pi^2}{2} \,(3 \,(\frac{\gamma}{\tau})^5\,+\, 4\, (\frac{\delta}{\tau})^5) \,=\,16.451... < 17.$$
\end{pf}

\subsection{Proof of Theorem~\ref{thm:weightedgaussian}}
Everything is now in place to complete the proof of Theorem~\ref{thm:weightedgaussian}. 
Assume $n \geq 5$, $d \geq 3$, and set $\epsilon \equiv n \bmod 2$ and $v={\rm v}_{n,d}$. 
By Corollary~\ref{cor:newfrmula}, we have 
$${\rm E}_{n,d}^H-\Theta^H_\epsilon(v) = {\rm T}^H_\epsilon(\widehat{\upmu_{n,d}}-\widehat{\gamma^\epsilon_v}).$$
Moreover, ${\rm T}^H_\epsilon : \mathcal{F} \rightarrow \C$ is linear and we have $\widehat{\upmu_{n,d}}-\widehat{\gamma^\epsilon_v}=A+B+C+D$ with $A,B,C,D \in \mathcal{F}$ as in 
Formula~\eqref{eq:decompfourier}, namely 
{\small 
$$A=  \frac{1}{{\rm u}_{n,d}} \bigl(\varphi_d(x)^n-e^{-vx^2/2}\bigr), \,\,\,B=\frac{1}{{\rm u}_{n,d}} \bigl(\Phi_{n,d}(x)-\Phi_{n,d}(0)\bigr),$$ 
$$C= -\frac{1}{{\rm M}_\epsilon(v)} \bigl({\rm G}^\star_v(x)-{\rm G}^\star_v(0)\bigr)\,\,\,{\rm and}\,\,\,D=\left(\frac{1}{{\rm u}_{n,d}}-\frac{1}{{\rm M}_\epsilon(v)}\right) \bigl(e^{-vx^2/2}-1\bigr).$$
 }
For $H\geq 0$, Lemma~\ref{lem:boundcentralterm} (i) and Prop.~\ref{prop:normTHeps} (ii) show ${\rm T}_\epsilon^H(A) \leq 0$. 
Moreover, we have $|{\rm T}_\epsilon^H(\phi)| \leq {\rm C}_\epsilon ||H||_\infty ||\phi||_\star$ for all 
$\phi \in \mathcal{F}$ by Prop.~\ref{prop:normTHeps} (i). Lemma~\ref{lem:boundcentralterm} (ii) shows $||A||_\star \leq {\rm M}_{n,d}$. Since $v\geq n/9\geq 5/9>4/\pi^2$, 
Lemmas~\ref{lem:boundPhi}, \ref{lem:boundGstar}, \ref{lem:lastterm} and \ref{lem:concattail} together show $||B||_\star+||C||_\star+||D||_\star \leq 17 n (7/8)^n$,  concluding the proof. $\square$

\section{Proof of Theorem~\ref{main:thm}}
\label{sect:pfmainthm}

\subsection{The functions $\Theta_\epsilon$ and $\Theta_0^\star$.}

In this section we assume once and for all that $H$ is the Odlyzko function ${\rm H}_{\log 2}$ of \S\ref{subsect:odlyzko}, and set $\Theta_\epsilon(v)=\Theta_\epsilon^{{\rm H}_{\log 2}}(v)$.

\begin{lemma} 
\label{lem:decroistheta}
For any $\epsilon \in \{0,1\}$, the map $v \mapsto \Theta_\epsilon(v)$ is differentiable and strictly decreasing on $]0,+\infty[$.
It tends to ${\rm J}(\epsilon)$ for $v\rightarrow 0$, and to $-\infty$ for $v \rightarrow +\infty$.
\end{lemma}

\begin{pf} 
For $v>0$, we have $\frac{\partial}{\partial v} e^{- \frac{i^2}{2v}} = \frac{i^2}{2v^2} e^{-\frac{ i^2}{2v}}$ and ${\rm m}_\epsilon(v)=\sum_{i \in \Z_\epsilon} e^{- \frac{i^2}{2v}}$, so  
$$2 v^2\, \frac{\partial}{\partial v}  \frac{e^{- \frac{j^2}{2v}} }{{\rm m}_\epsilon(v)} \,=\, 
\sum_{i \in \Z_\epsilon} j^2 \,\frac{e^{- \frac{j^2}{2v}} }{ {\rm m}_\epsilon(v)} \frac{e^{- \frac{i^2}{2v}} }{ {\rm m}_\epsilon(v)}\,-\, \sum_{i \in \Z_\epsilon} \,i^2 \,\frac{ e^{-\frac{i^2}{2v}}}{{\rm m}_\epsilon(v)} \frac{ e^{- \frac{j^2}{2v}}}{{\rm m}_\epsilon(v)}.$$
By multiplying by ${\rm J}(|j|)$ and summing over all $j \in \Z_\epsilon$ we deduce the equality $ 2 v^2 \Theta'_\epsilon(v) = \sum_{i,j \in \Z_\epsilon} (j^2-i^2) {\rm J}(|j|) \gamma_v^\epsilon(i) \gamma_v^\epsilon(j)$. Using the bijection $(i,j) \mapsto (j,i)$ of $\Z_\epsilon^2$ we finally obtain:
 $$ 4 v^2 \Theta'_\epsilon(v) \,=\, \sum_{i,j \in \Z_\epsilon} \bigl( j^2-i^2 \bigl) \,\bigl( {\rm J}(|j|) -{\rm J}(|i|) \bigl) \,\gamma_v^\epsilon(i) \,\gamma_v^\epsilon(j).$$
Each term in this sum is $\leq 0$ (and $<0$ for $|i|\neq |j|$) since ${\rm J}$ is strictly decreasing on $\R_{\geq 0}$, proving the first assertion. Note that the differentiations above are justified on every compact subinterval of $\R_{>0}$ by Lemma~\ref{lem:logcrois}.\par
 
For $v \rightarrow 0$, we clearly have ${\rm m}_0(v) \rightarrow 1$ and $\Theta_0(v) \rightarrow {\rm J}(0)$, 
as well as ${\rm  m}_1(v) \sim 2 e^{-1/2v}$ and $\Theta_1(v) \sim {\rm J}(1)$. Finally, for any fixed $N\geq 1$ we have 
$$\Theta_\epsilon(v) - {\rm J}(N) = \frac{1}{{\rm m}_\epsilon(v)} \sum_{j \in \Z_\epsilon} ({\rm J}(|j|)-{\rm J}(N)) e^{-j^2/2v}.$$
As ${\rm J}$ is decreasing, all the terms above with $|j|> N$ are $\leq 0$. 
Letting then $v \rightarrow +\infty$, and using ${\rm m}_\epsilon(v) \rightarrow +\infty$ from Lemma~\ref{lem:cargaussian}, we deduce
$${\rm limsup}_{v \rightarrow +\infty} \Theta_\epsilon(v)  \leq {\rm J}(N).$$
This being true for all $N\geq 1$, we conclude as ${\rm J}(x) \rightarrow -\infty$ for $x \rightarrow \infty$. \end{pf}

In the even case, we shall need the following variant of $\Theta_0(v)$, 
in which the index $i=0$ is omitted, namely:
$$\Theta_0^\star(v) = \sum_{i \in 2\Z \smallsetminus\{0\}} {\rm J}(|i|)\, \frac{e^{-\frac{i^2}{2v}}}{{\rm  m}_0^\star(v)},\,\,{\rm with}\,\, {\rm  m}_0^\star(v) = \sum_{i \in 2\Z \smallsetminus\{0\}}  e^{-\frac{i^2}{2v}} = {\rm  m}_0(v)-1.$$
The argument above applies verbatim,
with $(\Theta_0(v),{\rm  m}_0(v),\Z_0)$ replaced with 
$(\Theta_0^\star(v),{\rm  m}_0^\star(v),2\Z \smallsetminus\{0\})$, and shows:

\begin{lemma} 
\label{lem:decroisthetastar}
The map $v \mapsto \Theta_0^\star(v)$ is differentiable and strictly decreasing on $]0,+\infty[$. It tends to ${\rm J}(2)$ for $v\rightarrow 0$, and to $-\infty$ for $v \rightarrow +\infty$.
\end{lemma}

\begin{figure}[htbp]
\centering
\includegraphics[width=0.60\textwidth]{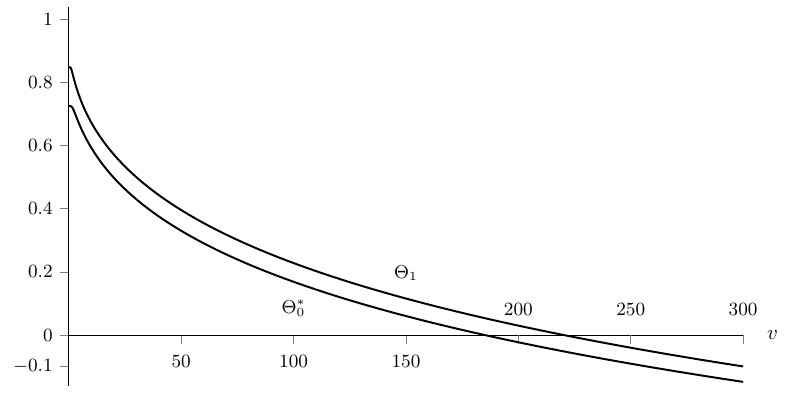}
\caption{The functions $\Theta_1$ and $\Theta_0^\star$.}
\label{fig:thetas}
\end{figure}

There is no difficulty to numerically compute $\Theta_\epsilon$ or $\Theta_0^\star$ : see Figure~\ref{fig:thetas} for their graphs on $]0,300[$. 
For the proof of the main theorem, we shall only need the following numerical values. Their verification is explained in \S\ref{subsect:numericalverif}.

\begin{lemma} 
\label{tab:valuestheta1}
The following values hold, rounded to $7$ decimal places:\ps\ps
\centering
\setlength{\tabcolsep}{4pt}
\begin{tabular}{c|ccc}
$3v$ & $650$ & $661$ & $661.34$ \\ \hline
$\Theta_1(v)$ & $0.0053942$ & $0.0001561$ & $-0.0000047$
\end{tabular}\hspace{1em}
\begin{tabular}{c|cc}
$3v$ & $550$ & $555$ \\ \hline 
$\Theta_0^\star(v)$ & $0.0030607$ & $0.0004149$
\end{tabular} 
\end{lemma}

We denote respectively by $\uprho_1$ and $\uprho_0^\star$ the unique zeros of $\Theta_1$ 
and $\Theta_0^\star$ on $]0,+\infty[$. 
These zeros exist by Lemmas~\ref{lem:decroistheta} and~\ref{lem:decroisthetastar}. 
A numerical computation shows 
$\uprho_1=220.4433355\ldots$ and $ \uprho_0^\star= 185.2624384\ldots$, and
$$3\uprho_1=661.3300065\ldots\, \, \, {\rm and}\, \, \,  3\uprho_0^\star=555.7873154\ldots.$$
As already emphasized in \S\ref{subsect:ideasandorg} for the odd case, 
they are very close to the constants in Theorem~\ref{main:thm} (i) and (ii) ! 
The main reason is the following version of Theorem~\ref{thm:weightedgaussian}.

\begin{lemma} 
\label{lem:explcompv90}
Let $n,d\geq 3$, $\epsilon \equiv n \bmod 2$ and $v = {\rm v}_{n,d}$. For $v\geq 90$ we have  
\begin{equation}
\label{eq:oddcomparison}
\Theta_\epsilon(v)-{\rm C}_\epsilon\,{\rm e}(v) -10^{-12} \ \leq\ {\rm E}_{n,d}\ \leq\ \Theta_\epsilon(v)+10^{-12},
\end{equation}
with ${\rm e}(v)= \frac{1}{15\,v}+\frac{64}{945\,v^2}$.
\end{lemma}

\begin{pf}  
We apply Theorem~\ref{thm:weightedgaussian}.
Note that we have ${\rm H}_{\log 2}\geq 0$ and $||{\rm H}_{\log 2}||_\infty=1$.
Moreover, using $x^4-1=(x+1)(x-1)(x^2+1)$ and $(d-1)^2+1 \leq d^2$ we find
$$ n \eta_d\,=\,n \,\frac{d(d-2)((d-1)^2+1)}{180d^4} \leq \frac{v}{60}.$$
Using $n \leq 9v$ we deduce ${\rm M}_{n,d} \leq {\rm e}(v)/{\rm u}_{n,d}$.
We also have $n > 3v \geq  270$, so
$$|{\rm u}_{n,d}^{-1}-1| \leq 1/((d-1)^{n-1}-1) \leq 1/ (2^{270}-1) < 2^{-269} < 10^{-26}.$$ 
Using ${\rm C}_\epsilon \leq 1$ and ${\rm e}(v) \leq {\rm e}(90) \leq 10^{-3}$ for $v\geq 90$, 
we deduce the inequality 
${\rm C}_\epsilon {\rm M}_{n,d} \,\leq \,{\rm C}_\epsilon {\rm e}(v) + 10^{-29}$ for all $n$. 
As $n (7/8)^n$ decreases with $n$ for $n \geq 7$, we have
${\rm R}_{n,d} \,\leq\, 17 \cdot 271\cdot (7/8)^{271} \,=\, 8.8634\ldots \,10^{-13}\, <\, 9 \cdot 10^{-13}$. 
\end{pf} 
  
\subsection{The odd case}

Let $n,d \geq 3$ with $n$ odd and $v={\rm v}_{n,d}$. Our aim is to prove Thm~\ref{main:thm} (i).\ps
	
We first claim ${\rm E}_{n,d}>0$ for $3v \leq 661$. Indeed, for $v\leq 100$, this is Proposition~\ref{prop:inegjensen} (i), so we may assume $v> 100$. By Lemma~\ref{lem:explcompv90} and ${\rm C}_1=1/2$, recall 
\begin{equation}
\label{eq:inegEodd90sub}
{\rm E}_{n,d}\ \geq \Theta_1(v)-\frac{1}{2}{\rm e}(v) -10^{-12}.
\end{equation}
We know that $ \Theta_1$ is decreasing by Lemma~\ref{lem:decroistheta}. For $v>100$ we have
$\frac{1}{2}{\rm e}(v) < \frac{1}{2}{\rm e}(100) < 3.4 \cdot 10^{-3}$. But by Lemma~\ref{tab:valuestheta1}, 
we also have $\Theta_1(650/3) > 5 \cdot 10^{-3} > 3.4 \cdot 10^{-3} + 10^{-12}$. 
We deduce ${\rm E}_{n,d}>0$ for $3v \leq 650$ by \eqref{eq:inegEodd90sub}. Assume finally $650 \leq 3v \leq 661$.
We have $\frac{1}{2} {\rm e}(650/3) < 1.55 \cdot 10^{-4}$ and $\Theta_1(661/3) > 1.56 \cdot 10^{-4} > 1.55 \cdot 10^{-4} + 10^{-12}$ by Lemma~\ref{tab:valuestheta1}, showing again ${\rm E}_{n,d} >0$ by \eqref{eq:inegEodd90sub}. \ps
We now claim ${\rm E}_{n,d}<0$ for $3v \geq 661.34$. Indeed, by Lemmas~\ref{lem:explcompv90} and~\ref{tab:valuestheta1} we have 
${\rm E}_{n,d}\ \leq \Theta_1(v)+10^{-12} \leq \Theta_1(661.34/3) + 10^{-12} < - 4 \cdot 10^{-6} + 10^{-12} < 0$.\ps
It only remains to determine the sign of ${\rm E}_{n,d}$ in the boundary strip 
$661<3v\leq 661.34$. We conclude the proof of Theorem~\ref{main:thm} (i) by the following lemma.
	
\begin{lemma} 
\label{lem:finitelistodd}
There are exactly $206$ pairs $(n,d)$ with $661 < n(d-2)/d \leq 661.34$ and $n$ odd. 
For those pairs, we have ${\rm E}_{n,d} > 0$ if, and only if, $3v<661.168$.
\end{lemma}

\begin{pf} For $(n,d)$ as in the statement we have $661< n \leq 1983$ and $d \leq 2n/(n-661.34)$ by Remark~\ref{rem:possvaluesv}, so those $(n,d)$ are straightforward to enumerate: we find $206$. The last ingredient is a simple numerical computation of ${\rm E}_{n,d}$ for those finitely many values,
for which we use the method described in \S\ref{subsect:numericalverif} below. 
The smallest value of $|{\rm E}_{n,d}|$ found is about $1.68\cdot 10^{-7}>10^{-10}$ (and is attained at $(n,d)=(663,724)$, so $3v = 661.16850...$). The $92$ pairs $(n,d)$ with the smallest values of $3v$ have a positive ${\rm E}_{n,d}$, and the $114$ remaining ones have a negative ${\rm E}_{n,d}$. The largest $3v$ among the positive is $\frac{229425}{347}=661.1671...$, attained at $(n,d)=(665,347)$, and the smallest $3v$ among the negative is $\frac{239343}{362}=661.1685...$, attained at $(n,d)=(663,724)$. As $661.168$ lies strictly between these two values, this concludes the proof of the lemma.
\end{pf}

\begin{remark} 
\label{rem:grh}
{\rm (On {\rm GRH} and the shape of the boundary)}
{\it
\begin{itemize}
\item[(i)] Assuming {\rm GRH} for $\Lambda(s,Z)$ in addition
to ${\rm HW}(Z)$ in the statement of Prop.~\ref{prop:maininequEF}, 
the condition {\rm (A3)} on $H$ may be weakened to
$\widehat H(\xi)\geq 0$ for all $\xi\in\R$. This is satisfied for 
the simpler Odlyzko function $\widetilde{{\rm H}}_\lambda(t)={\rm g}(t/\lambda)$ 
{\rm \cite[\S\S 2--3]{poitou}}. 
Using $\widetilde{{\rm H}}_\lambda \geq {\rm H}_\lambda$ and $\widetilde{{\rm H}}_\lambda(0)={\rm H}_\lambda(0)=1$,
we obtain ${\rm J}_{{\rm H}_\lambda} \leq {\rm J}_{\widetilde{{\rm H}}_\lambda}$, 
improving the resulting bounds.
In this context, this (standard) improvement is however quite modest.
Indeed, with the test function $\widetilde{{\rm H}}_{\log 2}$, 
computations give $666.0799$ instead of $661.168$ for Theorem~\ref{thmintro:HW} in the odd case, 
and $559.91$ instead of $555.6408$ in the even case {\rm (}see {\rm \cite[\S 4.3]{chetai}} 
for the numerical computation of ${\rm J}_{\widetilde{{\rm H}}_{\log 2}}${\rm )}.

\item[(ii)] During the computations above, for $n$ odd, we observed that there is no constant $C$
with the property that ${\rm E}_{n,d}^{\widetilde{{\rm H}}_{\log 2}}>0$ if and only if $n(d-2)/d < C$.  Indeed, a numerical computation gives
$${\rm E}_{667,1450}^{\widetilde{{\rm H}}_{\log 2}}\simeq -8.0385\cdot 10^{-8}  \,\,\, {\rm and} \,\,\, {\rm E}_{677,124}^{\widetilde{{\rm H}}_{\log 2}} \simeq 7.3543\cdot 10^{-7},$$
whereas the corresponding values of $n(d-2)/d$ are
respectively $666.08$ and $666.080645\ldots$.
This shows that the simple positivity criterion of
Theorem~\ref{main:thm} (i), especially for $n(d-2)/d$ close to $3\uprho_1$, is a specific feature 
of the weight ${\rm H}_{\rm log 2}$. 
\end{itemize}
}
\end{remark}

\subsection{The even case}

Our aim now is to prove Thm~\ref{main:thm} (ii). We will rely on the following variant of Lemma~\ref{lem:explcompv90}. We set ${\rm E}_{n,d}^\star={\rm E}_{n,d} - \upmu_{n,d}(0){\rm J}(0)$.
 
 \begin{lemma}
 \label{lem:evencompv90}
 Let $n,d\geq 3$ and $v={\rm v}_{n,d}$, with $n$ even. For $90 \leq v \leq 555/3$ we have ${\rm E}_{n,d}^\star  \geq \bigl( 1 - \sqrt{\frac{2}{\pi v}}\bigr) \Theta_0^\star(v) - \frac{\pi}{4} {\rm e}(v) - 10^{-11}$.
 \end{lemma}
 
 \begin{pf}
As a starting point, we have $ {\rm E}_{n,d} \geq \Theta_0(v) - \frac{\pi}{4} {\rm e}(v) - 10^{-12}$ by Lemma~\ref{lem:explcompv90}. By definition, we have ${\rm  m}_0^\star(v)/{\rm  m}_0(v)=1-1/{\rm  m}_0(v)$ hence
$$\Theta_0(v) \,=\, (1- 1/{\rm  m}_0(v))\,\Theta_0^\star(v)+\frac{{\rm J}(0)}{{\rm  m}_0(v)}.$$
Using ${\rm  m}_0(v) = \sqrt{\frac{\pi v}{2}} {\rm  M}_0(v)$ and ${\rm  M}_0(v) \geq 1$ by Lemma~\ref{lem:cargaussian} for $\epsilon=0$, note the inequality ${\rm  m}_0(v) \geq \sqrt{\frac{\pi v}{2}}$. As we have $\Theta_0^\star(v) \geq \Theta_0^\star(555/3)>0$ by Lemma~\ref{lem:decroisthetastar} and~\ref{tab:valuestheta1}, we deduce 
$$\Theta_0(v) \geq \bigl( 1 - \sqrt{\frac{2}{\pi v}}\bigr) \Theta_0^\star(v) \,+\, {\rm J}(0) \sqrt{\frac{2}{\pi v}} \,\bigl(\frac{1}{{\rm  M}_0(v)} - 1\bigr)\,+\, {\rm J}(0) \sqrt{\frac{2}{\pi v}}.$$
The last term is $\geq {\rm J}(0) \upmu_{n,d}(0)$ by Lemma~\ref{lem:boundmu} and ${\rm J}(0) >0$. 
So it only remains to check $|{\rm J}(0) \sqrt{\frac{2}{\pi v}} \,\bigl(\frac{1}{{\rm  M}_0(v)} - 1\bigr)|<10^{-12}$. It is actually much smaller: we have $$0 \leq 1-1/{\rm  M}_0(v) \leq {\rm  M}_0(v) -1 \leq 2 e^{-\pi^2v/2}/(1-e^{-\pi^2v/2}) < 10^{-190}$$ for $v\geq 90$ by the same argument as after Formula~\eqref{eq:majMepsv}, whereas $|{\rm J}(0)| \leq 1$ and $\sqrt{\frac{2}{\pi v}} \leq 1$.
 \end{pf}
 
 	We can now argue in a rather similar way as in the odd case.
Assume $n,d\geq 3$ with $n$ even, $(n,d) \neq (4,3)$ and set $v={\rm v}_{n,d}$. 
For $3v\leq 277$, we have ${\rm E}_{n,d}^\star >0$ by Proposition~\ref{prop:inegjensen} (ii), so we may assume $277 < 3v < 555.6408$. \ps
We assume first $277 < 3v\leq 555$ (so $v>90$). This ensures $
\Theta_0^\star(v)>0$ by Lemmas~\ref{lem:decroisthetastar} and~\ref{tab:valuestheta1}. Recall also that by the same lemma, $\Theta_0^\star$ is decreasing, whereas $v \mapsto 1 - \sqrt{\frac{2}{\pi v}}$ is increasing. 
We have $\frac{\pi}{4} {\rm e}(v) <\frac{\pi}{4} {\rm e}(90) < 6 \cdot 10^{-4}$. 
On the other hand, for $3v < 550$ we also have
$$\bigl( 1 - \sqrt{\frac{2}{\pi v}}\bigr) \Theta_0^\star(v) \geq (1-\sqrt{\frac{6}{277\pi}})  \Theta_0^\star(550/3) > 2.8 \cdot 10^{-3}.$$
This shows ${\rm E}_{n,d}^\star > 2 \cdot 10^{-3} > 0$ for $277 < 3v < 550$.\ps
Assume now $550 \leq 3v \leq 555$. We have then $\frac{\pi}{4} {\rm e}(v) \leq \frac{\pi}{4} {\rm e}(550/3) < 2.9 \cdot 10^{-4}$ and 
$$\bigl( 1 - \sqrt{\frac{2}{\pi v}}\bigr) \Theta_0^\star(v) \geq (1-\sqrt{\frac{6}{550 \pi}}) \Theta_0^\star(555/3) > 3.9 \cdot 10^{-4},$$
showing again ${\rm E}_{n,d}^\star > 9 \cdot 10^{-5} > 0$. \ps
It only remains to consider the pairs $(n,d)$ with $n$ even and 
$555<\frac{n(d-2)}{d}<555.6408$. As in Lemma~\ref{lem:finitelistodd}, 
these finitely many pairs are dealt with by a direct numerical computation.
	
\begin{lemma} 
\label{lem:finitelisteven}
There are exactly $2179$ pairs $(n,d)$ with $555 < n(d-2)/d < 555.6408$ and $n$ even,
and for all those pairs we have ${\rm E}^\star_{n,d} > 0$.
\end{lemma}

\begin{pf} Using the method described in \S\ref{subsect:numericalverif}, the smallest value found for ${\rm E}^\star_{n,d}$ is about $3.90\cdot 10^{-8}>10^{-10}$ and is attained at $(n,d)=(556,3095)$, so for $3v=555.64071\ldots$. The next possible value for $3v$ is $555.64082\ldots$ for $(n,d)=(556,3096)$, justifying the bound $555.6408$.
\end{pf}

\begin{remark} 
\label{rem:evenbalance}
By mimicking the proof in the odd case, we could also show that for $n$ even, we have\footnote{So the situation is as in Remark~\ref{rem:grh} (ii). An extra feature here is that the two pairs $(n,d)=(690,66)$ and $(n,d)=(736,22)$  
have the same $n(d-2)/d = 669.0909...$, but corresponding ${\rm E}_{n,d}$ with opposite signs.
For these computations, we can use the three values $\Theta_0(v)=0.00917..., \, \, 0.00058...$ and $-0.0000019...$, for $3v=650, 668$ and $669.25$ respectively: the remaining $1461$ pairs $(n,d)$ with $668< n(d-2)/d < 669.25$ and $n$ even (so $n \neq 669$) are treated by direct computation.
} 
${\rm E}_{n,d}>0 \Leftrightarrow n(d-2)/d < 669.0828$ or $(n,d)=(736,22)$. 
In particular, {\it if $n(d-2)/d < 669.0828$ then there is no hypersurface $Z$ of degree $d$ in $\mathbb{P}^{n-1}_{\Z}$ that is smooth over $\Z$, satisfying ${\rm HW}(Z)$, and that is
 balanced in the sense of Remark~\ref{rem:algcaseneven}. }
 \end{remark}

\subsection{Numerical verifications}
\label{subsect:numericalverif}

We now give the details of the numerical computations used in Lemmas~\ref{tab:valuestheta1}, \ref{lem:finitelistodd} and~\ref{lem:finitelisteven}. 
\medskip
\begin{center}
{\sc Computation of $\Theta_1$ and $\Theta_0^\star$}
\end{center}

Let us prove Lemma~\ref{tab:valuestheta1}.
We first explain how to compute $\Theta_1(v)$ with an error $<10^{-10}$ for $10 \leq v \leq 300$. 
It is easy to check $|{\rm J}'(x)| \leq \frac{{\rm log}\,2}{2}<1$ for $x\geq 1$, hence $|{\rm J}(x)| \leq x+|{\rm J}(1)|-1<x$ in this range. 
 As the function $t \mapsto t e^{-t^2/2v}$ is decreasing for $t\geq \sqrt{v}$, 
for each $N\geq 2$ and $v \leq (N-1)^2$ we have
$$ |\sum_{|i| \geq N} |{\rm J}(|i|)| e^{-\frac{i^2}{2v}} | \,< 2 \sum_{ i \geq N} i e^{-\frac{i^2}{2v}}   < \,2\,\int_{N-1}^\infty t\,e^{- \frac{t^2}{2v}} dt \,=\, 2 v\,e^{-\frac{(N-1)^2}{2v}}$$
In our range $10 \leq v \leq 300$, the right-hand side is $<10^{-13}$ for $N=150$. 
On the other hand, setting $\lambda(v) = \sqrt{\pi v/2}$ we have ${\rm  m}_1(v) = \lambda(v) {\rm M}_1(v)$ with $$|{\rm  M}_1(v) -1| \leq 2 e^{-\frac{v \pi^2}{2}}/(1- e^{-\frac{v \pi^2}{2}})$$ by Lemma~\ref{lem:cargaussian}, a quantity $< 10^{-21}$ for $v\geq 10$. In particular, we have ${\rm  m}_1(v) \geq 1$.
Set $\epsilon=10^{-11}$ and assume we have numbers $J_i$ with $|J_i - {\rm J}(|i|)| \leq \epsilon$ for each $1 \leq i < N$. For each odd $i$ in this range we have
$$| {\rm J}(|i|) \,\gamma_v^1(i) \,-\, J_i \,\frac{e^{-\frac{i^2}{2v}}}{\lambda(v)} |  \,=\, \gamma_v^1(i) \,| {\rm J}(|i|)\,-\,J_i \,{\rm  M}_1(v)|$$ 
$$\leq \gamma_v^1(i) \,( \epsilon \,+\, |J_i|\, |{\rm  M}_1(v)-1| )\, \leq \,\gamma_v^1(i)\,( \epsilon\, +\, (N-1+\epsilon) \,|{\rm  M}_1(v)-1|).$$
Using finally $\sum_{1\leq |i| < N, \, \, i \,\,{\rm odd}} \gamma_v^1(i) \leq 1$, we deduce the finite approximation
$$|\Theta_1(v)\, -\, \frac{2}{\lambda(v)}\, \sum_{1 \leq i < 150, i\,{\rm odd}} \,J_{i}\, e^{-\frac{i^2}{2v}} | < 10^{-13} + 150 \cdot 10^{-21} + 10^{-11} < 10^{-10}.$$ 
The exact same inequality holds with $\Theta_1$, ``odd'' and $\,\lambda(v)$ replaced with 
$\Theta_0^\star$, ``even'' and $\,\lambda(v)-1$. The argument is the same: setting ${\rm  M}_0^\star(v) = {\rm  m}_0^\star(v)/(\lambda(v)-1)$, we have 
$|{\rm  M}_0^\star(v)-1| = \frac{\lambda(v)}{\lambda(v)-1} |{\rm  M}_0(v)-1| < 2 \frac{\lambda(v)}{\lambda(v)-1}  e^{-v \pi^2/2}/(1- e^{-v \pi^2/2})$, which is again $<10^{-21}$ for $v\geq 10$. For the tail, note that by Lemma~\ref{lem:cargaussian} we still have ${\rm  m}_0^\star(v) =  \lambda(v) {\rm  M}_0(v) -1  \geq \lambda(v)-1>1$ for $v\geq 10$. $\square$
\medskip
\begin{center}
{\sc Computation of ${\rm E}_{n,d}$ and
  ${\rm E}_{n,d}^\star$}
\end{center}

We finally explain how we computed ${\rm E}_{n,d}$ and
${\rm E}_{n,d}^\star$ in Lemmas~\ref{lem:finitelistodd} 
and~\ref{lem:finitelisteven}. Our (non-optimized) method consists of computing
the Hodge numbers exactly, and of approximating the masses
and the values of ${\rm J}$ with controlled errors:

\begin{itemize}
\item[(i)] Note that ${\rm h}_{n,d}(i)=0$ for $|i|> n(d-2)/d$, since we clearly have 
$|S_{n,d}| \leq n(d-2)/d$, and we may also assume $i\leq 0$ by the symmetry ${\rm h}_{n,d}(i)={\rm h}_{n,d}(-i)$. For $n,d$ fixed, the computation of the integers ${\rm h}_{n,d}(i)$ with $|i| \leq n(d-2)/d$ is done by induction. Indeed, 
if $c(k)$ denotes the coefficient of $t^k$ in $Q(t)=(1+t+t^2+\dots+t^{d-2})^n$, we have 
${\rm h}_{n,d}(i)=c((n+i)d/2 -n)$ for all $i \in n+2\Z$ by definition. Taking the logarithmic derivative 
of $Q(t)=(1-t^{d-1})^n(1-t)^{-n}$ leads to the equality
$(1-t-t^{d-1}+t^d)Q'(t)=n(1-(d-1)t^{d-2}+(d-2)t^{d-1})Q(t)$,
and comparing the coefficients of $t^{k-1}$ gives $k c(k) = (k+n-1)c(k-1) + (k-(n+1)(d-1))c(k-d+1)+ (n(d-2)+d-k)c(k-d)$ for $k\geq 1$. It thus is enough to determine $c(k)$ for $k \leq dn/2 - n$.\ps
\item[(ii)] For our purposes we have $n(d-2)/d<662$ in all cases, so the relevant values of ${\rm J}$ are the ${\rm J}(i)$ with $0 \leq i \leq 661$. We mention the upper bound $\sum_{0 \leq |i| \leq 661} |{\rm J}(|i|)| < 4000$.
\ps
\item[(iii)] The ${\rm h}_{n,d}(i)$, like $\upbeta_{n,d}$, can be large integers, with several thousands of digits in our range, so it is better in the calculation to compute the rational number $\upmu_{n,d}(i) = {\rm h}_{n,d}(i)/\upbeta_{n,d}$ up to a certain precision. Given $\epsilon,\delta>0$, assume we have real numbers $\mu_i$ and $J_i$ for each $i$ in the support $S$ of $\upmu_{n,d}$, with $|\mu_i-\upmu_{n,d}(i)| \leq  \delta$ and $|J_i-{\rm J}(|i|)| \leq \epsilon$. Then we have 
$$|{\rm E}_{n,d}-\sum_{i \in S} \mu_i J_i| \leq \epsilon + \delta \epsilon |S| + \delta \sum_{i \in S} |{\rm J}(|i|)|,$$
using $\upmu_{n,d}(i) {\rm J}(|i|)-\mu_i J_i=\upmu_{n,d}(i)({\rm J}(|i|)-J_i)+(\upmu_{n,d}(i)-\mu_i)(J_i-{\rm J}(|i|))+(\upmu_{n,d}(i)-\mu_i){\rm J}(|i|)$ and $\sum_{i \in S} \upmu_{n,d}(i)=1$.
\end{itemize}

For the purpose of computing ${\rm E}_{n,d}$ with $n(d-2)/d <662$, and with an error $<10^{-10}$, 
it is thus enough to take $\epsilon=10^{-11}$ and $\delta=10^{-14}$. The computation of $\sum_{i \in S} \mu_i J_i$ is exact, and can even be done using only integers if, say, the $\mu_i$ and $J_i$ are given as decimal numbers. For $n$ even, the same method computes ${\rm E}_{n,d}^\star$ with the same error bound, simply by omitting the index $i=0$ in the sums above, since we still have $\sum_{i \in S, i \neq 0} \upmu_{n,d}(i) \leq 1$.

\section{Proof of Theorem~\ref{thm:autexist}}
\label{sect:persppf}

In this section we prove Theorem~\ref{thm:autexist}.
Let $\N_{\rm odd} \subset \N$ denote the subset of odd integers. 
For any nonempty finite subset $F \subset \N_{\rm odd}$, 
we denote by ${\rm S}(F)$ the set of level $1$, algebraic regular, selfdual symplectic, cuspidal automorphic representations of $\GL_{2|F|}$ over $\Q$ with weights $\{ \pm i/2\,\,|\,\,i \in F\}$, and set ${\rm N}(F) = |{\rm S}(F)|$. The following lemma is essentially tautological.

\begin{lemma} 
\label{lem:comblem}
Assume $n,d \geq 3$ with $n$ odd and $f \geq 1$. Then we have $\Pi_{n,d,f} \neq \emptyset$ if, and only if, there is a family $\mathcal{F}$ of nonempty subsets $F \subset \mathbb{N}_{\rm odd}$, together with a collection of positive integers $({\rm m}_F)_{F \in \mathcal{F}}$, such that:
\begin{itemize}
\item[(i)] For all $F \in \mathcal{F}$, we have $|F| \leq f$ and ${\rm N}(F)  \neq 0$, 
\item[(ii)] For all $i \in \N_{\rm odd}$, $\sum_{F \in \mathcal{F}\, |\, i \in F} {\rm m}_F\, =\, {\rm h}_{n,d}(i)$.
\end{itemize}
\end{lemma}

\begin{pf} 
Assume we have $\mathcal{F}$ and $({\rm m}_F)_{F \in \mathcal{F}}$ satisfying (i) and (ii).
Then we may choose arbitrarily some $\pi_F \in {\rm S}(F)$ by (i), and set $\pi_\mathcal{F}=\boxplus_{F \in \mathcal{F}} \pi_F^{\boxplus {\rm m}_F}$. Then $\pi_\mathcal{F}$ trivially belongs to $\Pi_{n,d,f}$ by (i) and (ii).
Conversely, assume we have $\pi=\boxplus_{j=1}^s \pi_j$ in $\Pi_{n,d,f}$. For $1 \leq j \leq s$, define $F_j$ as the set of $i \in \N_{\rm odd}$ such that $i/2$ is a weight of $\pi_j$. Set $\mathcal{F}=\{ F_j\, |\, 1\leq j \leq s\}$, and for $F \in \mathcal{F}$, define ${\rm m}_F$ as the number of $1 \leq j \leq s$ with $F=F_j$. Then (i) and (ii) hold now by definition of $\Pi_{n,d,f}$. 
\end{pf}

We now recall a few known facts about the integers ${\rm N}(F)$.
If we have $F=\{a_1,\dots,a_k\}$ with $a_1 < \dots < a_k$ we also set 
${\rm N}(a_1,\dots,a_k)={\rm N}(F)$. \ps

Of course, for $a \in \N_{\rm odd}$, ${\rm N}(a)$ is the dimension of the space of modular cuspforms of weight $a+1$ for ${\rm SL}_2(\Z)$. By~\cite[Prop. 4.1]{chrenard}, we also have 
\begin{equation}
\label{eq:formSab} {\rm N}(a,b) \,=\, {\rm d}(a,b) \,-\, \delta_{a=1\,\, \&\,\, b \equiv 1 \bmod 4} \,{\rm N}(b),
\end{equation}
where ${\rm d}(a,b)$ denotes the dimension of 
the space of Siegel modular cuspforms for ${\rm Sp}_4(\Z)$ 
of weight ${\rm Sym}^{a-1}\otimes\det^{(b-a)/2+2}$, 
the subtraction above corresponding to the contribution of Saito-Kurokawa lifts. 
Classical works of Igusa for $a=1$ and of Tsushima for $a>1$ and $b-a\geq 6$, 
completed by Ta\"ibi and Petersen independently for $b-a \leq 4$, 
give an explicit formula for the dimension ${\rm d}(a,b)$: 
see {\it e.g.} \cite[\S 9]{vandergeer} and \cite{tsushima, petersen, taibidimtrace}. \ps
More generally, an explicit formula for the map $F \mapsto {\rm N}(F)$ 
has been obtained in \cite{chrenard,taibidimtrace,chetai} for $|F| \leq 8$. 
Beyond the special cases $|F| \leq 2$ recalled above, 
the only other cases of these formulas that we shall need below are $3 \leq |F| \leq 6$, 
for which we refer to \cite{chrenard,taibidimtrace}. 
We emphasize that the proofs of these formulas, 
as well as that of \eqref{eq:formSab}, 
rely on Arthur's endoscopic classification \cite{arthur}, 
which has been conditional for some time, 
but whose missing parts have been established 
in the recent series of works \cite{MW,AGIKMS,HD}.\ps

\subsection{Proof of assertion (i) of Theorem~\ref{thm:autexist}.} 

Set $A = \{1,3,5,7,9,13\}$. All we shall need about the formulas for ${\rm N}(a)$ and ${\rm N}(a,b)$ is contained in the following lemma.

\begin{lemma} 
\label{lem:dimavailable}
Let $1\leq a<b$ be odd integers. \begin{itemize}
\item[(i)] We have ${\rm N}(a) \neq 0 \iff a \notin A$.\par
\item[(ii)] For $a\in A$, we have ${\rm N}(a,b) \neq 0 \iff$ $b$ is as in Table~\ref{tab:Sab}.\par
\item[(iii)] We have ${\rm N}(a,b)= 0$ for $b \in A$.\ps
\end{itemize}
\vspace{-.3cm}
\begin{table}[htbp]
\centering
  \small
  \setlength{\tabcolsep}{5pt}
  \begin{tabular}{c|cc}
  $a$ & $b\equiv 1 \bmod 4$ & $b\equiv 3 \bmod 4$ \\
  \hline
  $1$ & $b\geq 37$ & $b=67$ or $b\geq 75$ \\
  $3$ & $b\geq 41$ & $b\geq 27$ \\
  $5$ & $b\geq 21$ & $b\geq 31$
  \end{tabular}
  \qquad
  \begin{tabular}{c|cc}
  $a$ & $b\equiv 1 \bmod 4$ & $b\equiv 3 \bmod 4$ \\
  \hline
  $7$  & $b\geq 25$ & $b\geq 19$ \\
  $9$  & $b\geq 21$ & $b\geq 23$ \\
  $13$ & $b\geq 21$ & $b\geq 23$
  \end{tabular}
\caption{The odd integers $b$ for which ${\rm N}(a,b)\neq 0$, for $a \in A$.}
\label{tab:Sab}
  \end{table}
  \end{lemma}
  
\begin{pf} Assertion (i) is clear. There are many ways to verify assertion (ii). 
One is to observe that we have ${\rm d}(a,b) \leq {\rm d}(a,b+8)$ by considering the multiplication by the weight $4$ genus $2$ Eisenstein series ${\rm E}_4$. A numerical application of the formulas above for ${\rm d}(a,b)$ allows us to verify Table~\ref{tab:Sab} for all $b \leq 82$ and $a \in A$, and (iii). For $a>1$ or $b\equiv 3 \bmod 4$, it shows ${\rm N}(a,b)\neq 0$ for $75\leq b \leq 82$, hence ${\rm N}(a,b) \neq 0$ for all such $b\geq 75$. So it only remains to consider 
$a=1$ and $b \equiv 1 \bmod 4$, that is $k:=(b-1)/2+2$ even. For $k$ even Igusa's formula 
implies ${\rm d}(1,b) \geq |\mathcal{Q}_k|$ with 
$$\mathcal{Q}_k = \{ (u,v,s,t) \in \mathbb{N}^4\, \, | \,\,4u+6v+10s+12t = k\,\,{\rm  and} \,\,s+t \geq 1\}.$$
(See~\cite[Thm. 9.1]{vandergeer}.) But if ${\rm M}_q({\rm SL}_2(\Z))$ denotes the space of modular forms of weight $q$ for ${\rm SL}_2(\Z)$, whose dimension is well-known to be $|\{(u,v) \in \N^2 |\, 4u+6v=q\}|$, 
we have ${\rm N}(b) = \dim {\rm M}_{2k-2-12} = \dim {\rm M}_{k-10} + \dim {\rm M}_{k-12}$.
So 
${\rm N}(b)$ coincides with the number of $(u,v,s,t)\in \mathcal{Q}_k$ with $s+t=1$. 
But for even $k\geq 20$, {\it i.e.}
$b \geq 37$, there is $(u,v,s,t) \in \mathcal{Q}_k$ with $v=0$ and either $(s,t)=(2,0)$ or $(s,t)=(1,1)$, 
showing ${\rm N}(1,b)>0$.
\end{pf}

For each $a \in A$ we denote by ${\rm B}_a$ the set of odd integers $b>a$ with ${\rm N}(a,b) \neq 0$, {\it i.e.} given by Table~\ref{tab:Sab}. Note that for all $a \in A$ and $b \in {\rm B}_a$ we have $a \leq 13$ and $b\geq 19$, so $A \cap {\rm B}_a =\emptyset$. 

\begin{lemma}
\label{lem:hallmarriage}
Let $(h_i)_{i \in \N_{\rm odd}}$ be a nonzero sequence of nonnegative integers with finite support. 
We have an equivalence between:
\begin{itemize}
\item[(i)] there is a family $\mathcal{F}$ of nonempty subsets $F \subset \mathbb{N}_{\rm odd}$, together with a collection of positive integers $({\rm m}_F)_{F \in \mathcal{F}}$, such that: for all $F \in \mathcal{F}$, we have $|F| \leq 2$ and ${\rm N}(F)  \neq 0$, and for all $i \in \N_{\rm odd}$, we have $\sum_{F \in \mathcal{F}\, | i \in F}{\rm m}_F= h_i$. \ps
\item[(ii)] For all $I \subset A$, we have $\sum_{i \in I} h_i \leq \sum_{i \in {\rm B}_I} h_i$, with ${\rm B}_I=\cup_{a \in I} {\rm B}_a$.
\end{itemize}
\end{lemma}

\begin{pf} 
This is a special case of Hall's marriage theorem~\cite{hall}. 
Indeed, consider a finite set $S$ and a map $h: S \rightarrow \N_{\rm odd}$ 
satisfying $|h^{-1}(i)|=h_i$ for all $i \in \N_{\rm odd}$.
Set $S_A = S \cap h^{-1}(A)$.
By Lemma~\ref{lem:dimavailable}, assertion (i) is equivalent to 
the existence of an injective map $f: S_A \rightarrow S$, 
such that $h(f(s)) \in {\rm B}_{h(s)}$ for all $s \in S_A$.
By Hall's marriage theorem, such a marriage exists if, and only if, 
for each nonempty subset $X \subset S_A$ 
we have $|X| \leq |h^{-1}({\rm B}_{h(X)})|$. 
In such an inequality, the right-hand side only depends on $I:=h(X)$ 
and coincides with $\sum_{i \in {\rm B}_I} h_i$, 
whereas the left-hand side is maximal when $X=h^{-1}(I)$, 
in which case it coincides with $\sum_{i \in I} h_i$. 
This shows (i) $\Leftrightarrow$ (ii).
\end{pf}

We finally assume $d=3$. By definition (see the proof of Lemma~\ref{lem:griffitsproba}), 
recall that ${\rm h}_{n,3}(i)$ is the coefficient of $u^{3i}$ in the Laurent 
polynomial $(u+u^{-1})^n$, {\it i.e}
$${\rm h}_{n,3}(i) = \binom{n}{(n+3i)/2}$$
for all $i \in \Z$, with the convention that the right-hand side is $0$ when 
$(n+3i)/2$ is not an integer in $[0,n]$. Note also $\upbeta_{n,3}=(2^n-2)/3$ by Formula~\eqref{eq:formbetand}.
Assertion (i) of Theorem~\ref{thm:autexist} follows then from Lemmas~\ref{lem:hallmarriage} \& \ref{lem:comblem}, and from the following lemma.

\begin{lemma} 
\label{lem:calccubic}
Assume $n\geq 3$ is odd and set $h_i={\rm h}_{n,3}(i)$ for $i \in \mathbb{N}_{\rm odd}$.
Then:\begin{itemize}
\item[(i)] Condition (ii) of Lemma~\ref{lem:hallmarriage} holds for all $n\geq 4893$.
\item[(ii)] We have $h_1 > \sum_{i \in {\rm B}_1} h_i$ for all $n\leq 4891$.
\end{itemize}
\end{lemma}

Note that due to the Gaussian shape of the ${\rm h}_{n,3}(i)$, it was expected that (i) should happen for $n$ large enough. 

\begin{pf} 
It is straightforward to check with the computer that the $2^{|A|}=64$ inequalities in condition (ii) of Lemma~\ref{lem:hallmarriage} hold for $n=n_0=4893$.\footnote{The inequalities even hold if we replace each ${\rm B}_I$ with ${\rm B}_I \cap [19,99]$. The integers ${\rm h}_{n,3}(i)$ being 
rather large (but less than $(2^n-2)/3$), we can save memory by computing instead an approximation of the rational ${\rm h}_{n,3}(i)/\upbeta_{n,3}$: to justify the inequalities, it is enough to compute them with an error $\leq 10^{-8}$. } We claim that these inequalities hold then a fortiori for all $n\geq n_0$. Indeed, assuming ${\rm h}_{n,3}(i) \neq 0$ we have 
$${\rm q}_n(i):=\frac{{\rm h}_{n+2,3}(i)}{{\rm h}_{n,3}(i)} \,=\, 4 \frac{(n+1)(n+2)}{(n+2)^2-9i^2},$$
since both $n+3i$ and $2n-(n+3i)=n-3i$ increase by $2$ when $n$ increases by $2$.
So $i \mapsto {\rm q}_n(i)$ is increasing on $\N_{\rm odd} \cap [0,n/3]$. As a consequence, for any $\emptyset \neq I \subset A$, the sum $\sum_{i \in I} {\rm h}_{n,3}(i)$ is at most multiplied by ${\rm q}_n(13)$ when $n$ increases by $2$ and $n/3 \geq 13$, since ${\rm max}\, A=13$. 
Similarly, when $n$ increases by $2$ the sum $\sum_{i \in {\rm B}_I} {\rm h}_{n,3}(i)$ is at least multiplied by ${\rm q}_n(19)$ since ${\rm Min}\, {\rm B}_I \geq 19$ and adding new elements $i$ in the support of ${\rm h}_{n+2,3}(i)$ can only increase the sum. Assertion (i) follows then by induction using ${\rm q}_n(13)<{\rm q}_n(19)$ for $n \geq n_0$ (since $n_0/3 \geq 19$). \par
Similarly, the inequality in assertion (ii) is obvious for $n < 3 \cdot 37 = 111$ (since ${\rm Min} \,{\rm B}_1=37$) and follows from a straightforward computer check for $n=4891$. It still holds by descending induction for all $111 \leq n \leq 4891$ by the same argument as above: when 
$n$ decreases by $2$, ${\rm h}_{n,3}(1)$ is divided by ${\rm q}_{n-2}(1)$ whereas $\sum_{i \in {\rm B}_1} {\rm h}_{n,3}(i)$ is divided at least by ${\rm q}_{n-2}(37)>{\rm q}_{n-2}(1)$.
\end{pf}

\subsection{Proof of assertion (ii).}

In order to prove this assertion, we shall use the following criterion.

\begin{lemma} 
\label{lem:weakhall}
Let $(h_i)_{i \in \N_{\rm odd}}$ be a nonzero sequence of nonnegative integers with finite support. Set $C=\sum_i h_i$ and $m= {\rm Max}_i h_i$. If $C \geq 27 m$ then there is a family $\mathcal{F}$ and collection $({\rm m}_F)_{F \in \mathcal{F}}$ satisfying assertion (i) of Lemma~\ref{lem:hallmarriage}.
\end{lemma}

\begin{pf} 
For $I \subset A$, we set ${\rm B}^{\rm c}_I = \N_{\rm odd} \smallsetminus {\rm B}_I$.
An inspection of Table~\ref{tab:Sab} shows $|{\rm B}^{\rm c}_{\{1\}}|=26$, 
and $|{\rm B}^{\rm c}_{\{a\}}| \leq |\N_{\rm odd}\cap [1,41]| \leq 21$ for $a \in A \smallsetminus\{1\}$. 
That being said, fix a nonempty $I \subset A$. 
We clearly have $\sum_{i \in I} h_i \,\leq \,|I|\,m$ and
$\sum_{i \in {\rm B}_I} h_i\, \geq \,C - |{\rm B}^{\rm c}_I|m$. 
It is thus enough to show $|I|\,m\, \leq \,C \,-\, |{\rm B}^{\rm c}_I|\,m$, or more strongly,
\begin{equation}
\label{eq:inegCmtoshow} 
|I| \,+ \,| {\rm B}^{\rm c}_I| \,\,\leq \,\,27
\end{equation}
using $27 m \leq C$. If $I=\{1\}$, 
we have $|I|+| {\rm B}^{\rm c}_I|=1+26$, 
so the equality holds in \eqref{eq:inegCmtoshow}.
Otherwise, $I$ contains some $a>1$, so we have 
$|{\rm B}^{\rm c}_I| \leq |{\rm B}^{\rm c}_{\{a\}}| \leq 21$, 
and $|I| \leq |A|=6$, showing $|I| + |{\rm B}^{\rm c}_I| \leq 6+21=27$ again.
\end{pf}

We now prove the first part of assertion (ii) of Theorem~\ref{thm:autexist}.
We may assume $d\geq 4$ by assertion (i).
Set $h_i={\rm h}_{n,d}(i)$.
We have $\sum_{i \in \N_{\rm odd}} h_i = \upbeta_{n,d}/2$, so
Lemma~\ref{lem:maxupmu} below shows that the assumptions of~\ref{lem:weakhall} are satisfied.
We conclude by Lemmas~\ref{lem:comblem} and \ref{lem:weakhall}.

\begin{lemma} 
\label{lem:maxupmu}
For all $n\geq 11139$ odd, $d\geq 4$, and $i \in \Z$, we have 
$\upmu_{n,d}(i) \leq 1/54$. 
\end{lemma}

\begin{pf} By Lemma~\ref{lem:boundmu}, we have
$\upmu_{n,d}(i) \leq \frac{1}{{\rm u}_{n,d}} \sqrt{\frac{2}{\pi {\rm v}_{n,d}}} \leq \frac{1}{1-1/3^{n-1}} \sqrt{\frac{12}{\pi n}}$, 
using ${\rm u}_{n,d} \geq 1-1/3^{n-1}$ for $n$ odd and ${\rm v}_{n,d} \geq n/6$ for $d\geq 4$.
We conclude as $n\mapsto \frac{1}{1-1/3^{n-1}} \sqrt{\frac{12}{\pi n}}$ is decreasing, and is $0.018517\ldots < 1/54$ for $n=11139$. 
\end{pf}

In order to prove the second part of assertion (ii), we start with two preliminary lemmas.

\begin{lemma} 
\label{lem:comblem2}
Assume we have $n,d \geq 3$ with $n$ odd and $f=2$. Suppose that there is $\mathcal{F}$ and $({\rm m}_F)_{F \in \mathcal{F}}$ satisfying the assumptions (i) and (ii) of Lemma~\ref{lem:comblem}, as well as 
${\rm N}(F)\, \, |\, \, {\rm m}_F$ for all $F \in \mathcal{F}$. Then $\Pi_{n,d,2}^\Q \neq \emptyset$.
\end{lemma}

\begin{pf}
Set $\pi_\mathcal{F} = \boxplus_{F \in \mathcal{F}} \Pi_F^{\,\boxplus\,{\rm m}_F/{\rm N}(F)}$, where $\Pi_F = \boxplus_{\varpi \in {\rm S}(F)} \varpi$. The fact that for any $F \subset \N_{\rm odd}$ with $|F|\leq 2$, the Satake parameters of $\Pi_F|\cdot |^{w/2}$, for any $w \in \Z$ odd, are defined over $\Q$, 
follows from the standard fact that the spaces of modular cuspforms for ${\rm SL}_2(\Z)$, and of Siegel modular cuspforms for ${\rm Sp}_4(\Z)$, have a natural $\Q$-structure preserved by the Hecke operators (the Saito–Kurokawa subspace being defined over $\Q$ as well). The necessity of the odd $|.|^{1/2}$ twist is a standard feature of the Satake isomorphism.
\end{pf}
\begin{lemma} 
\label{lem:divhndi}
Fix $n\geq 3$ odd and an integer $N\geq 1$. 
There are infinitely many $d\geq 3$ such that for all $i \in \Z$ we have ${\rm h}_{n,d}(i) \equiv 0 \bmod N$.
\end{lemma}

\begin{pf} 
It is well-known that for each odd integer $i \in ]-n,n[$, there is a polynomial ${\rm P}_{n,i} \in \Q[X]$ such that
${\rm h}_{n,d}(i)={\rm P}_{n,i}(d)$ for all $d\geq 3$ \cite[Cor. 17.3.3]{arapura}. We claim $(n-1)! \,{\rm P}_{n,i}(X) \in (X-2)\Z[X]$. This claim implies ${\rm h}_{n,d}(i) \equiv 0 \bmod N$ for $d \equiv 2 \bmod N(n-1)!$, hence the lemma. \ps
Set $p+1=(n+i)/2$. Griffiths' formula recalled in \S\ref{subsect:hodgehypersurfaces} asserts that ${\rm h}_{n,d}(i)$ coincides with the coefficient ${\rm c}_{n,p}(d)$ in $t^{(p+1)d-n}$ 
of $(1-t^{d-1})^n(1-t)^{-n} \in \Z[t]$. Using $(1-t)^{-n} = \sum_{r \geq 0} \binom{r+n-1}{n-1} t^r$, $(1-t^{d-1})^n=\sum_{s=0}^n (-1)^s \binom{n}{s} t^{(d-1)s}$, we find that for all integers $d\geq 1$, we have
{\small $${\rm c}_{n,p}(d) = \sum_{s=0}^{p}(-1)^s\binom{n}{s}
  \binom{r(s,d)+n-1}{n-1},$$}
\par \noindent with $r(s,d):=(p+1-s)d+s-n$. Indeed, for $0 \leq s \leq n$ and $r\geq 0$, the equality $r+(d-1)s=(p+1)d-n$ is equivalent to $r=(p+1-s)(d-1)+p+1-n$, hence to $r=r(s,d)$. Note that $r\geq 0$ forces $s\leq p$, since $p+1<n$.
Conversely, for $0 \leq s \leq p$ and $d\geq 1$, we have $r(s,d)+n-1 \geq 0$, so the binomial $\binom{r(s,d)+n-1}{n-1}$ vanishes if $r(s,d)< 0$, establishing the formula above for ${\rm c}_{n,p}(d)$. This formula is clearly a polynomial in $d$ with coefficients in $\frac{1}{(n-1)!}\Z$, reproving the existence of ${\rm P}_{n,i}$ with $(n-1)!\, {\rm P}_{n,i} \in \Z[X]$. For $d=2$, we have 
$(1-t^{d-1})^n(1-t)^{-n}=1$ and $(p+1)d-n$ odd, showing ${\rm c}_{n,p}(2)={\rm P}_{n,i}(2)=0$.
\end{pf}

We finally prove the second part of assertion (ii) of Theorem~\ref{thm:autexist}.
Fix $n\geq 11139$.
Choose an integer $N\geq 1$ which is divisible by all the nonzero values of ${\rm N}(a)$ for $1 \leq a < n$ odd,
and by all the nonzero values of ${\rm N}(a,b)$ for $1\leq a<b<n$ odd.
By Lemma~\ref{lem:divhndi}, there are infinitely many $d\geq 3$ satisfying ${\rm h}_{n,d}(i) \equiv 0 \bmod N$ for all $i \in \Z$. 
Choose such a $d \geq 4$. 
In particular, we have $\upbeta_{n,d} \equiv 0 \bmod N$. 
For $i\geq 1$ odd, we set $h_i={\rm h}_{n,d}(i)/N$. 
The quantities $C=\sum_{i} h_i=\frac{\upbeta_{n,d}}{2N}$ and 
$m\,=\,{\rm Max}_i \,h_i = ({\rm Max}_i\,{\rm h}_{n,d}(i))/N$ still satisfy 
$$m/C \,=\,2\, {\rm Max}_i\,\upmu_{n,d}(i) \leq 1/27$$
by Lemma~\ref{lem:maxupmu}. 
By applying Lemma~\ref{lem:weakhall} to this sequence $(h_i)$,
we obtain a family $\mathcal{F}$ and a collection $({\rm m}_F)_{F \in \mathcal{F}}$ 
as in assertion (i) of Lemma~\ref{lem:hallmarriage}. 
We conclude by Lemma~\ref{lem:comblem2} applied to 
the family $\mathcal{F}$ and the collection $(N {\rm m}_F)_{F \in \mathcal{F}}$. 
$\square$

\subsection{Proof of assertion (iii)}

All we need is to provide a family $\mathcal{F}$ and 
a collection of integers $({\rm m}_F)_{F \in \mathcal{F}}$ as in Lemma~\ref{lem:comblem}, 
in the special case $n=2773$ (so $n/3=924+1/3$), $d=3$ and $f=6$.
We simply describe such a family.
The $43$ elements $F \in \mathcal{F}$ with $|F|>1$ are listed in Table~\ref{tab:profiles2773} below.
The remaining $419$ members of $\mathcal{F}$ are the singletons $\{i\}$ for $i=11$ and for all odd $89\leq i\leq 923$.
For the singletons $\{i\}$ other than $\{89\}$, we simply set ${\rm m}_{\{i\}}={\rm h}_{n,3}(i)=\binom{n}{(n+3i)/2}$.
The other multiplicities ${\rm m}_F$ satisfy the system of linear equations
$$ \sum_{ F\in\mathcal{F}\, |\, i\in F}\, {\rm m}_F
  \,=\,{\rm h}_{n,3}(i),\,\,\,{\rm with}\,\, 1 \leq i \leq 89\, \,\,{\rm odd}\,\,{\rm and} \, \, i \neq 11.$$
This system has thus $44$ equations (the number of $i$ above), 
$44$ unknowns (the $F\in \mathcal{F}$ with $|F|>1$ or $F=\{89\}$), and turns out 
to be invertible with determinant $\pm 24$.  
An exact computation shows that its unique solution consists of positive integers.
Note that those integers are all insanely huge: they all have between $827$ and $833$ digits. 
This is why we give approximate values of ${\rm m}_F/\upbeta_{n,3}$ instead in Table~\ref{tab:profiles2773}.
The values of ${\rm N}(F)$ are computed using the formulae
recalled above for $1\leq |F|\leq 2$, and those of~\cite{taibidimtrace} for $3 \leq |F| \leq 6$.
\begin{table}[htbp]
\centering
\small
\setlength{\tabcolsep}{1.5pt}
\renewcommand{\arraystretch}{1.08}
\resizebox{\linewidth}{!}{%
\begin{tabular}{@{}c|c|r@{\hspace{3pt}\vrule width 0.8pt\hspace{3pt}}c|c|r@{\hspace{3pt}\vrule width 0.8pt\hspace{3pt}}c|c|r@{\hspace{3pt}\vrule width 0.8pt\hspace{3pt}}c|c|r@{}}
$F$ & ${\rm N}(F)$ & ${\rm m}_F/\upbeta_{n,3}$
& $F$ & ${\rm N}(F)$ & ${\rm m}_F/\upbeta_{n,3}$
& $F$ & ${\rm N}(F)$ & ${\rm m}_F/\upbeta_{n,3}$
& $F$ & ${\rm N}(F)$ & ${\rm m}_F/\upbeta_{n,3}$ \\
\hline
$\{89\}$ & 7 & 0.00000009
& $\{1,7,81\}$ & 17 & 0.00000107
& $\{1,7,9,63\}$ & 1 & 0.00007230
& $\{3,5,9,17,41\}$ & 1 & 0.00174875 \\
$\{7,19\}$ & 1 & 0.01574628
& $\{1,5,83\}$ & 4 & 0.00000063
& $\{1,3,15,65\}$ & 1 & 0.00004770
& $\{1,5,13,17,43\}$ & 2 & 0.00226167 \\
$\{5,21\}$ & 1 & 0.02222436
& $\{5,9,17,25\}$ & 1 & 0.00015755
& $\{1,7,9,67\}$ & 1 & 0.00003106
& $\{1,7,9,17,45\}$ & 1 & 0.00169957 \\
$\{1,85\}$ & 14 & 0.00000036
& $\{3,9,17,27\}$ & 1 & 0.00655052
& $\{1,3,15,69\}$ & 2 & 0.00001996
& $\{1,7,9,17,49\}$ & 3 & 0.00092303 \\
$\{1,87\}$ & 2 & 0.00000021
& $\{3,9,15,29\}$ & 1 & 0.00017144
& $\{1,7,9,71\}$ & 3 & 0.00001266
& $\{1,7,9,15,51\}$ & 1 & 0.00066707 \\
$\{1,89\}$ & 16 & 0.00000003
& $\{3,7,17,29\}$ & 1 & 0.00793808
& $\{3,5,7,73\}$ & 1 & 0.00000793
& $\{1,3,13,17,53\}$ & 1 & 0.00047586 \\
$\{3,15,23\}$ & 1 & 0.00756005
& $\{1,9,15,35\}$ & 1 & 0.00622745
& $\{1,7,9,75\}$ & 5 & 0.00000490
& $\{1,5,9,17,55\}$ & 2 & 0.00033507 \\
$\{3,13,25\}$ & 1 & 0.01633059
& $\{1,9,13,41\}$ & 1 & 0.00074550
& $\{1,5,9,77\}$ & 1 & 0.00000299
& $\{1,5,9,17,59\}$ & 7 & 0.00015976 \\
$\{5,9,27\}$ & 1 & 0.00737729
& $\{1,7,15,41\}$ & 1 & 0.00047652
& $\{1,7,13,19,31\}$ & 1 & 0.00955831
& $\{1,5,9,15,23,33\}$ & 1 & 0.00776546 \\
$\{1,15,29\}$ & 1 & 0.00350371
& $\{1,5,13,57\}$ & 1 & 0.00023288
& $\{1,7,13,17,37\}$ & 2 & 0.00098780
& $\{1,3,9,17,23,37\}$ & 1 & 0.00394176 \\
$\{1,5,79\}$ & 3 & 0.00000180
& $\{1,5,13,61\}$ & 3 & 0.00010818
& $\{1,7,13,15,39\}$ & 1 & 0.00385178
& $\{1,5,9,15,17,47\}$ & 2 & 0.00126067 \\
\end{tabular}%
}
\caption{\small The $44$ elements $F \in \mathcal{F}$ with $|F|>1$ or $F=\{89\}$, their associated
${\rm N}(F)$ and normalized multiplicities ${\rm m}_F/\upbeta_{n,3}$
for $n=2773$ (rounded to eight decimal places).}
\label{tab:profiles2773}
\end{table}

\begin{remark} 
For $f=3,4,5$, 
we also found similar examples in $\Pi_{n,3,f}$ with $n=3069$, $2823$ and $2781$ respectively, {\it i.e.} with $n/3=1023$, $941$ and $927$. 
\end{remark}

\end{document}